\documentclass[10pt,a4paper]{article}%
\usepackage[centertags]{amsmath}
\usepackage{amsfonts}
\usepackage{amssymb}
\usepackage{amsthm}
\usepackage{epsfig}
\usepackage{setspace}
\usepackage{ae}
\usepackage{eucal}
\usepackage[usenames]{color}%
\usepackage{graphicx}

\theoremstyle{plain}
\newtheorem{thm}{Theorem}[section]
\newtheorem{lem}[thm]{Lemma}
\newtheorem{cor}[thm]{Corollary}
\newtheorem{prop}[thm]{Proposition}
\theoremstyle{definition}
\newtheorem{defi}[thm]{Definition}

\theoremstyle{remark}

\newtheorem{rem}[thm]{Remark}

\begin{document}

\title{Weak singularities of 3-D Euler equations and restricted regularity of Navier Stokes equation solutions with time dependent force terms }
\author{J\"org Kampen }
\maketitle

\begin{abstract}
A class of singular solutions of the $3$D Euler equation corresponds to a class of singular solutions of the incompressible Navier Stokes equation with time dependent force terms, as the time-dependent force terms can be chosen such that they cancel the viscosity terms of the incompressible Navier Stokes equation. 
In this context classical solution branches of the three dimensional incompressible Euler equation are constructed where at least one vorticity component blows up at some point after finite time for regular velocity component data. Furthermore, we show that there are classical solution branches with $C^k\cap H^k$-data for $k\geq 2$ which develop weak singularities (kinks) at some point of space time in $C^{k-1}\setminus C^{k}$ for any $k\geq 2$ after finite time. The time-local solution branches of the Euler equation are obtained by viscosity  limits of solutions of Navier-Stokes-type viscosity ($\nu$)-extensions of time-reversed Euler-type equations.  Here, Lipschitz continuity of (spatial first order derivatives of) the Leray projection term and the Burgers term are used in local time iteration schemes in terms of convolutions with the Gaussian of dispersion $\nu$ and first order spatial derivatives of this Gaussian.  For the time reversed equations and related Navier Stokes type equation extensions data with singular vorticity at the origin and strong polynomial decay are used. The functional solution increment and derivatives up to some order, i.e., the solution minus the spatially convoluted data snd spatial derivatives of this increment up to some order, preserve strong polynomial decay of some order. This polynomial decay of the local time solution is strong enough such that compactness arguments can be constructed in classical Banach spaces. The Lipschitz continuity of the first order derivatives of the Leray projection term and of the Burgers term implies that regularity is gained after any finite time for the time-reversed equation. This result can be sharpened in the sense that the second order derivatives of the Leray projection term of a local solution are Lipschitz even in a viscosity limit. As a consequence, the original Euler equation has a solution branch which develops singular solutions or kinks from regular data after some finite time.

\end{abstract}


\section{Idea of construction and statement of weak singularity theorems }
We determine data with weakly singular short-time solutions of the three dimensional incompressible Euler equation on the whole Euclidean space. Let $v_i,~1\leq i\leq D$ denote the velocity component functions in dimension $D=3$, and let
\begin{equation}
\omega =\mbox{curl}(v)=\left(\frac{\partial v_3}{\partial x_2}-\frac{\partial v_2}{\partial x_3},\frac{\partial v_1}{\partial x_3}-\frac{\partial v_3}{\partial x_1},\frac{\partial v_2}{\partial x_1}-\frac{\partial v_1}{\partial x_2} \right)
\end{equation}
denote the corresponding vorticity functions. The incompressible Euler equation in vorticity form is
\begin{equation}\label{vorticity}
\frac{\partial \omega}{\partial \tau}+v\cdot \nabla \omega=\frac{1}{2}\left(\nabla v+\nabla v^T\right)\omega, 
\end{equation}
and for the corresponding Cauchy problem this equation has to be solved with some initial data $\mbox{curl}(f)$ at time $\tau=0$ of some function $f=(f_1,f_2,f_3)^T$ with $f_i\in {\mathbb R}^3\rightarrow {\mathbb R}$. It is well-known that (for sufficiently regular $\omega$) the velocity is determined by
\begin{equation}\label{vel}
v(\tau,x)=\int_{{\mathbb R}^3}K_3(x-y)\omega(\tau,y)dy,~\mbox{where}~K_3(x)h=\frac{1}{4\pi}\frac{x\times h}{|x|^3}.
\end{equation}
This Biot-Savart law is also important for our construction of a viscosity limit of an extended equation of Navier Stokes type. This limit is not a unique but a local time solution of the Euler vorticity equation which is Lipschitz and of some order of polynomial decay at spatial infinity. The Biot-Savart law can then be imposed in order to obtain Lipschitz continuity (and some order of polynomial decay at spatial infinity) of first order spatial derivatives of the velocity components for some positive time, which in turn leads to higher regularity of the velocity solution increment (velocity solution minus initial data) in a viscosity limit.     
In order to construct local-time and long-time singularities we consider first a related equation which is obtained from (\ref{vorticity}) by the simple time transformation
\begin{equation}
\tau=-t,~\omega^-_i(t,.)=\omega_i(\tau,.),~v^-_i(t,.)=v_i(\tau,.),~1\leq i\leq D.
\end{equation}
Spatial derivatives are untouched, so multiplying the spatial part of transformed equation by $-1$ and writing the equation in coordinates we have
\begin{equation}\label{vort1}
 \frac{\partial \omega^-_i}{\partial t}-\sum_{j=1}^3v^-_j\frac{\partial \omega^-_i}{\partial x_j}=-\sum_{j=1}^3\frac{1}{2}\left(\frac{\partial v^-_i}{\partial x_j}+\frac{\partial v^-_j}{\partial x_i} \right)\omega^-_j,~1\leq i\leq D=3.
\end{equation}
We refer to the equation in (\ref{vort1}) as the 'time-reversed Euler equation'.
We shall observe that the 'Navier-Stokes type' extension of the latter equation with an additional viscosity term, i.e., the equation 
\begin{equation}\label{vort2}
 \frac{\partial \omega^{\nu,-}_i}{\partial t}-\nu\Delta \omega^{\nu,-}_i-\sum_{j=1}^3v^-_j\frac{\partial \omega^{\nu,-}_i}{\partial x_j}=-\sum_{j=1}^3\frac{1}{2}\left(\frac{\partial v^{\nu,-}_i}{\partial x_j}+\frac{\partial v^{\nu,-}_j}{\partial x_i} \right)\omega^{\nu,-}_j,
\end{equation}
for $1\leq i\leq D$ (and with a positive constant $\nu >0$), has a regular solution on a time interval $[0,T]$ for some $T>0$ for a certain class of data in $H^1$ which have a singular vorticity at some point in space time with spatial dimension $n\geq 3$. Regularity results for this specific class of data $\omega^{\nu,-}_i(0,.)=\omega^{f,-}_i,~1\leq i\leq D$ lead to the conclusion that short or long-time solution families with parameter $\nu>0$ have a viscosity limit solution branch 
\begin{equation}
\omega^{-}_i(t,.)=\lim_{\nu\downarrow 0}\omega^{\nu,-}_i(t,.),~1\leq i\leq D
\end{equation}
of the time-reversed Euler equation with the same initial data, where some regularity is preserved. The initial data $\omega^{f,-}_i,~1\leq i\leq D$ of the the time-reversed Euler equation are then 'final data' of the original incompressible Euler equation, which, hence, has a vorticity singularity after finite time. If the solution branch constructed for the time-reversed Euler equation is short-time, then the solution for the original Euler equation is short-time with a short-time singularity, and if the solution branch constructed for the time-reversed Euler equation is long-time, then the solution for the original Euler equation is long-time with a long-time singularity. We state these results more precisely in the two theorems at the end of this section. In this paper we shall prove the existence of short time singularities for regular data. The proof of long-time singularities follows from the construction of global regular upper bounds, which imply the existence of global regular solution branches. This extension to global time is not the main interest of this paper, so we sketch an argument for the existence of global regular solution branches of the Euler equation in a remark below. In the following let $C^{\infty}$ be the function space of real-valued functions on ${\mathbb R}^D$ with continuous derivatives of any order, and let $H^s$ denote the Sobolev function space (in $L^2$-sense) of order $s\geq 0$ (for real numbers $s$). Furthermore, we denote the space of functions (defined on the domain $[0,T)\times {\mathbb R}^D$) with continuous time derivatives up to order $m$ (nonnegative integer) and multivariate spatial derivatives up to order $n$ (nonnegative integer) by $C^{m,n}\left(\left[0,T\right)\times {\mathbb R}^D \right)$ or just by $C^{m,n}$ if the reference to the domain is known from the context. We have
\begin{thm}\label{main1}
Let $D=3$. There exist data $\omega_i^f=\omega_i(0,.)\in H^2\cap C^{2},~1\leq i\leq D$ and a vorticity solution $\omega_i,~1\leq i\leq D$ of the three dimensional incompressible Euler equation Cauchy problem such that after some finite time $T>0$ there is a blow-up of the classical solution, i.e.,
\begin{itemize}
 \item[i)] for some $T>0$ there is a solution function $\omega_i:[0,T)\times {\mathbb R}^D\rightarrow {\mathbb R},~1\leq i\leq D$ in $C^{1}\left(\left[0,T\right), H^2\cap C^2 \right)$ which satisfies the incompressible Euler equation pointwise on the domain $\left[0,T\right)\times {\mathbb R}^D$ in a classical sense with data $\omega_i(0,.)=\omega_i^f,~1\leq i\leq D$;
 \item[ii)] for the solution in item i) we have
 \begin{equation}\label{supomega}
 \sup_{\tau\in [0,T)}|\omega_i(\tau,x)|= \infty,
 \end{equation}
i.e., there is no finite upper bound for the left side of (\ref{supomega})on the time interval $[0,T)$. 
\end{itemize}
\end{thm}
Even a long time singularity version of Theorem \ref{main1} holds. We have
\begin{thm}\label{main2}
A stronger version of Theorem \ref{main1} holds with the same text except for the replacement of the quantor 'after some finite time $T>0$' by 'after any finite time $T>0$'.
\end{thm}
\begin{rem}\label{globrem}
Theorem \ref{main2} follows from Theorem \ref{main1} if the existence of global regular solution branches of the Euler equation can be shown. Strictly speaking following our argument we need to construct global regular solution branches of the time-reversed Euler equation for regular data. However, the following argument does not depend on the signs of the nonlinear terms such that it suffices to construct global regular solution branches of the incompressible Euler equation. In the next remark we note some simple observations which lead to the construction of global regular solution branches.
\end{rem}
\begin{rem}
It is not a simple scaling argument which leads to the conclusion that a certain classes of nonlinear equations have global regular solution branches. We sketch an argument here which may be applied for a certain class of equations which have no damping term. In order to construct global regular solution branches we have to exploit some structure of the equations.  A first essential characteristic of the classes of equations we have in mind is that local time solutions can be represented in terms of convolutions with first order spatial derivatives of the Gaussian. In case of the incompressible Navier Stokes equation the incompressibility condition helps in order to obtain local time representations of the value function itself in terms of convolutions with the first order spatial derivatives of the Gaussian (note that this is not possible for typical blow-up equations as $v_{i,\tau}=v_{,i}^2,~ 1\leq i\leq D$).  For the construction of global regular solution branches of the incompressible Euler equation we need the construction of a viscosity limit $\nu \downarrow 0$ in addition. This is possible if the convoluted data of the nonlinear terms with the Gaussian and or with first order spatial derivatives of the Gaussian in a local time representation of the solution of the Navier Stokes equation are Lipschitz. This leads to upper bounds in terms of second moments of the Gaussian, where the upper bound of this second moment constant is an essential difference to other equations which do not have local representations in terms of convolutions with spatial first order derivatives of the Gaussian. We note that the dimension constraint $D\geq 3$ is due to structural constraints in the case of doemension less than three which are apparent in the the vorticity for ofthe equation.
Let us now sketch an argument starting with the local time solution of the Navier Stokes equation. We introduce the scaling $\tau=\rho \tau_{\rho}$ which is also useful later when we consider local contraction results and their independence from the  viscosity parameter $\nu >0$. Let $v^{\rho}_i(\tau_{\rho},.)=v_i(\tau,.)$. We have 
\begin{equation}
v_{i,\tau_{\rho}}=\rho v_{i,\tau}.
\end{equation}
Let us assume that for some time $\tau_{\rho 0}\geq 0$ we have
\begin{equation}\label{vrhoinit}
v^{\rho}_i(\tau_{\rho 0},.)\in H^m\cap C^m,~{\big |}v^{\rho}_i(\tau_{\rho 0},.){\big |}_{H^m\cap C^m}\leq C_m ~\mbox{ for some }~ m\geq 2.
\end{equation}
For the construction of global regular solution branches the parameter $\rho >0$ will depend only on the spatial ${\big |}.{\big |}_{H^m\cap C^m}$-norm and on a constant related to the local $L^1$-norms of the /first order derivatives of the) Laplacian kernel, i.e., on
\begin{equation}
C_K:=1+\max_{1\leq i\leq D}\int_{B_1(0)}{\big |}K_{,i}(y){\big |}dy+\int_{B_1(0)}{\big |}K(y){\big |}dy
+{\big |}K_{D,i}{\big |}_{L^2({\mathbb R}^D\setminus B_1)},
\end{equation}
where $B_1(0)$ is the ball of radius $1$ around the origin in ${\mathbb R}^D$.
Next we note that local contraction arguments lead to the conclusion of existence of local time solutions in a time interval $[\tau_0,\tau_0+\Delta]$ for some $\Delta >0$.
For $1\leq i\leq D$, and multiindices 
$\beta,\gamma$ with $1\leq |\beta|=|\gamma|+1$, where 
$\beta_k=\gamma_k+1$ and
 $\beta_j=\gamma_j$ for $j\in \lbrace 1,\cdots,D\rbrace\setminus \lbrace k\rbrace$ we have for all $1\leq i\leq D$ the local time representation
\begin{equation}\label{Navlerayscheme3rca}
\begin{array}{ll}
 D^{\beta}_xv^{\rho,\nu}_i(\tau_{\rho},.)=D^{\beta}_xv^{\rho,\nu}_i(t_0,.)\ast_{sp}G^{\rho}_{\nu}(\tau_{\rho},.)\\
\\
- \rho\int_{0}^{\tau_{\rho}}\int_{{\mathbb R}^D}\sum_{j=1}^D D^{\gamma}_x\left( v^{\rho,\nu}_jv^{\rho,\nu}_{i,j}\right) (\sigma,y)\left( \frac{-2(.-y)_k}{4 \nu\rho  (\tau_{\rho}-\sigma)}\right)  G^{\rho}_{\nu}dyd\sigma\\
\\+\rho\int_{0}^{\tau_{\rho}}\int_{{\mathbb R}^D}\left( \int_{{\mathbb R}^D}\left( K_{D,i}(y-z)\right) \sum_{j,m=1}^DD^{\gamma}_x\left( \frac{\partial v^{\rho,\nu}_m}{\partial x_j}\frac{\partial v^{\rho,\nu}_j}{\partial x_m}\right) (\sigma,z)dz\right)\times \\ 
\\
\times 
\left( \frac{-2(.-y)_k}{4\nu  \rho (\tau_{\rho}-\sigma)}\right) G^{\rho}_{\nu}(\tau_{\rho}-\sigma,.-y)dyd\sigma,
\end{array}
\end{equation}
where the first derivatives of the spatially scaled Gaussian $G^{\rho}_{\nu}$ are given by
\begin{equation}
\begin{array}{ll}
G^{\rho}_{\nu,k}(\tau_{\rho},z;0,0)=\left( \frac{-2z_k}{4\rho \nu\tau_{\rho}}\right) \frac{1}{\sqrt{4\pi \rho\nu \tau_{\rho}}^D}\exp\left(-\frac{|z|^2}{4\rho\nu \tau_{\rho}} \right)\\
\\
=:\left( \frac{-2z_k}{4\rho \nu\tau_{\rho}}\right)G^{\rho}_{\nu}(\tau_{\rho},z;0,0)=:\left( \frac{-2z_k}{4\rho \nu\tau_{\rho}}\right)G^{\rho}_{\nu}(\tau_{\rho},z),
\end{array}
\end{equation}
where the latter definition is an abbreviation sometimes used in the following.
Now it is interesting that the value functions themselves have a local time representation in terms of convolutions with the first order spatial derivatives of the Gaussian. Indeed, due to incompressibility we have
\begin{equation}\label{Navlerayscheme3rcb}
\begin{array}{ll}
 v^{\rho,\nu}_i(\tau_{\rho},.)=v^{\rho,\nu}_i(t_0,.)\ast_{sp}G^{\rho}_{\nu}(\tau_{\rho},.)\\
\\
- \rho\int_{0}^{\tau_{\rho}}\int_{{\mathbb R}^D}\sum_{j=1}^D \left( v^{\rho,\nu}_jv^{\rho,\nu}_i\right) (\sigma,y)\left( \frac{-2(.-y)_i}{4 \nu \rho(\tau_{\rho}-\sigma)}\right)  G^{\rho}_{\nu}dyd\sigma\\
\\+\rho\int_{0}^{\tau_{\rho}}\int_{{\mathbb R}^D}\left( \int_{{\mathbb R}^D}\left( K_D(y-z)\right) \sum_{j,m=1}^D\left( \frac{\partial v^{\rho,\nu}_m}{\partial x_j}\frac{\partial v^{\rho,\nu}_j}{\partial x_m}\right) (\sigma,z)dz\right)\times \\ 
\\
\times 
\left( \frac{-2(.-y)_i}{4\nu  \rho (\tau_{\rho}-\sigma)}\right) G^{\rho}_{\nu}(\tau_{\rho}-\sigma;.-y)dyd\sigma.
\end{array}
\end{equation} 
Here the incompressibility condition is used in order to rewrite the Burgers term. We may abbreviate
\begin{equation}
B^0_i\ast G^{\rho}_{\nu}:=\int_{0}^{\tau_{\rho}}\int_{{\mathbb R}^D}\sum_{j=1}^D \left( v^{\rho,\nu}_jv^{\rho,\nu}_i\right) (\sigma,y)\left( \frac{-2(.-y)_i}{4 \nu \rho(\tau_{\rho}-\sigma)}\right)  G^{\rho}_{\nu}dyd\sigma,
\end{equation}
and
\begin{equation}
\begin{array}{ll}
L^0\ast G^{\rho}_{\nu}:=\int_{0}^{\tau_{\rho}}\int_{{\mathbb R}^D}\left( \int_{{\mathbb R}^D}\left( K_D(y-z)\right) \sum_{j,m=1}^D\left( \frac{\partial v^{\rho,\nu}_m}{\partial x_j}\frac{\partial v^{\rho,\nu}_j}{\partial x_m}\right) (\sigma,z)dz\right)\times \\ 
\\
\times 
\left( \frac{-2(.-y)_i}{4\nu  \rho (\tau_{\rho}-\sigma)}\right) G^{\rho}_{\nu}(\tau_{\rho}-\sigma;.-y)dyd\sigma,
\end{array}
\end{equation}
which implicitly defines the symbols $B^0_i$ and $L^0$, where we note that the latter does not depend on the component index $i$. From the equation (\ref{Navlerayscheme3rca}) we may define abbreviations $B_i$ and $L_i$ for the convoluted Burgers and Leray data analogously. Spatial derivatives of these functions are denoted in Einstein notation.
We observe that the representation in (\ref{Navlerayscheme3rcb}) as well as the representation in (\ref{Navlerayscheme3rca}) have the parameter $\rho$ as a coefficient. Furthermore for local time regular solutions $v^{\rho,\nu}_i(\tau_{\rho},.)\in H^m\cap C^m,~ 1\leq i\leq D$ with $m\geq 2$ and $\tau_{\rho}\in [\tau_{\rho 0},\tau_{\rho 0}+\Delta]$ the Burgers-and Leray functions $B^0_i,L^0,~ 1\leq i\leq D$, $B_i,L_i,~ 1\leq i\leq D$, and $B_{i,j},L_{i,j},~ 1\leq i,j\leq D$ are all Lipschitz such that we get interesting upper bounds of the convolutions
\begin{equation}
B^0_i\ast G^{\rho}_{\nu,k}, L^0\ast G^{\rho}_{\nu,k},~B_i\ast G_{\nu,k}, L_i\ast G_{\nu,k},~ B_{i,j}\ast G^{\rho}_{\nu,k}, L_{i,j}\ast G^{\rho}_{\nu,k}.
\end{equation}
Note that the first order spatial derivatives of the Gaussian are antisymmetric in the sense that
\begin{equation}
\left( \frac{-2y_i}{4\nu  \rho (\tau_{\rho}-\sigma)}\right) G^{\rho}_{\nu}(\tau_{\rho},y)=\left( \frac{2y^-_i}{4\nu  \rho (\tau_{\rho}-\sigma)}\right) G^{\rho}_{\nu}(\tau_{\rho},y^-),
\end{equation}
where $y^-=(y^-_1,\cdots,y^-_D)$ with $y^-_i=-y_i$ and $y^-_j=y_j$ for $j\in \lbrace 1,\cdots D\rbrace\setminus \lbrace i\rbrace$.
For any spatially global Lipschitz continuous function $F$ we have
\begin{equation}\label{estest}
\begin{array}{ll}
{\Big |}\rho F\ast G^{\rho}_{\nu,i}=\rho\int F(\tau_{\rho-}\sigma,x-y)\left( \frac{-2y_i}{4\nu  \rho(\sigma)}\right) G^{\rho}_{\nu}(\sigma,y)dyd\sigma {\Big |}\\
\\
={\Big |}\rho \int \int_{y_i\geq 0} \left( F(\tau_{\rho}-\sigma,x-y)-F(\tau_{\rho}-\sigma,x-y^-)\right)\left( \frac{-2y_i}{4\nu  \rho (\tau_{\rho}-\sigma)}\right) G^{\rho}_{\nu}(\sigma,y)dyd\sigma{\Big |}\\
\\
\leq {\Big |}\rho L\int \int_{y_i\geq 0} \left( \frac{4y^2_i}{4\nu  \rho \sigma}\right) G^{\rho}_{\nu}(\sigma,y)dyd\sigma {\Big |}=4\rho LM_2,
\end{array}
\end{equation} 
where $M_2$ is a finite second moment constant of the Gaussian which is small in the following sense (note here, that the integral $\int_{y_i\geq 0}$ is with respect to the domain $\lbrace y\in {\mathbb R}^D|y_i\geq 0\rbrace$) . We consider the case $D=3$. First  note that for any small time $\sigma>0$ the main mass of the spatial integral is on a ball of radius $\sqrt{\rho \nu}$. Now, given $\nu >0$, small $\sqrt{\rho\nu}$ and small time $\sigma$, and using partial integration and polar coordinates we observe that (denoting the three-dimensional ball by $B^3_{\sqrt{\rho\nu}}$ and the one dimensional ball by $B_{\sqrt{\rho\nu}}=\left\lbrace r>0|r\leq \sqrt{\nu r}\right\rbrace $, and choosing implicitly a time interval $[\tau_{\rho 0},\tau_{\rho 0}+\Delta])$ with $\tau_{\rho 0}=0$ without loss of generality ) 
\begin{equation}\label{series}
\begin{array}{ll}
\frac{1}{4\pi^2}{\Big |}\rho L\int \int_{\left\lbrace y|  y_i\geq 0\right\rbrace \cap B^3_{\sqrt{\rho\nu}}} \left( \frac{4y^2_i}{4\nu  \rho \sigma}\right) G^{\rho}_{\nu}(\sigma,y)dyd\sigma {\Big |}\\
\\
\leq {\Big |}\rho L\int \int_{r\in B_{\sqrt{\rho\nu}}} \left( \frac{4r^2}{4\nu  \rho \sigma}\right) G^{\rho}_{\nu}(\sigma,r)r^2drd\sigma {\Big |}
\\
\\
\leq {\Big |}\rho L\int   \left( \frac{r^5}{5\nu  \rho \sigma}\right) \frac{1}{\sqrt{4\pi \rho\nu \sigma}^D}\exp\left(-\frac{|r|^2}{4\rho\nu \sigma} \right){\Big |}^{\sqrt{\rho\nu}}_{0}d\sigma \\
\\
+\rho L\int \int_{r\in  B_{\sqrt{\rho\nu}}} \left( \frac{r^5}{5\nu  \rho \sigma}\right) \left( \frac{-2r}{4\nu  \rho\sigma}\right)\frac{1}{\sqrt{4\pi \rho\nu \sigma}^D}\exp\left(-\frac{r^2}{4\rho\nu \sigma} \right)drd\sigma {\Big |}\\
\\
\leq {\Big |}\rho L\int   \frac{1}{5\sigma} \frac{1}{\sqrt{4\pi \sigma}^D}\exp\left(-\frac{1}{4\sigma} \right)d\sigma \\
\\
-\rho L\int \int_{r\in  B_{\sqrt{\rho\nu}}} \left( \frac{r^6}{10\nu^2  \rho^2 \sigma^2}\right) \frac{1}{\sqrt{4\pi \rho\nu \sigma}^D}\exp\left(-\frac{r^2}{4\rho\nu \sigma} \right)drd\sigma {\Big |}.
\end{array}
\end{equation}

Iterating this partial integration procedure we get a geometric series upper bound.
More precisely, note that for the last term on the right side in (\ref{series}) we have (recall $D=3$)
\begin{equation}\label{series2}
\begin{array}{ll}
 {\Big |}\rho L\int \int_{r\in  B_{\sqrt{\rho\nu}}} \left( \frac{r^6}{10\nu^2  \rho^2 \sigma^2}\right) \frac{1}{\sqrt{4\pi \rho\nu \sigma}^D}\exp\left(-\frac{r^2}{4\rho\nu \sigma} \right)drd\sigma {\Big |}\\
\\
\leq {\Big |}\rho L\int   \left( \frac{r^7}{70\nu^2 \rho^2 \sigma^2}\right) \frac{1}{\sqrt{4\pi \rho\nu \sigma}^D}\exp\left(-\frac{|r|^2}{4\rho\nu \sigma} \right){\Big |}^{\sqrt{\rho\nu}}_{0}d\sigma {\Big |}\\
\\
+\rho L\int \int_{r\in  B_{\sqrt{\rho\nu}}} \left( \frac{r^7}{70\nu^2  \rho^2 \sigma^2}\right) \left( \frac{-2r}{4\nu  \rho\sigma}\right)\frac{1}{\sqrt{4\pi \rho\nu \sigma}^D}\exp\left(-\frac{r^2}{4\rho\nu \sigma} \right)drd\sigma {\Big |}\\
\\
\leq {\Big |}\rho L\int   \frac{1}{70\sigma^2} \frac{1}{\sqrt{4\pi \sigma}^D}\exp\left(-\frac{1}{4\sigma} \right)d\sigma \\
\\
-\rho L\int \int_{r\in  B_{\sqrt{\rho\nu}}} \left( \frac{r^8}{140\nu^3  \rho^3 \sigma^3}\right) \frac{1}{\sqrt{4\pi \rho\nu \sigma}^D}\exp\left(-\frac{r^2}{4\rho\nu \sigma} \right)drd\sigma {\Big |}.
\end{array}
\end{equation}
Inductively we have for any $m\geq 1$
\begin{equation}\label{series3}
\begin{array}{ll}
{\Big |}\rho L\int \int_{\left\lbrace y|  y_i\geq 0\right\rbrace \cap B_{\sqrt{\rho\nu}}} \left( \frac{4y^2_i}{4\nu  \rho \sigma}\right) G^{\rho}_{\nu}(\sigma,y)dyd\sigma {\Big |}\\
\\
\leq {\Big |}\sum_{k=1}^m(-1)^{k+1}\rho L\int   \frac{R_k}{\sigma^k} \frac{1}{\sqrt{4\pi \sigma}^D}\exp\left(-\frac{1}{4\sigma} \right)d\sigma {\Big |}\\
\\
+ {\Big |}\rho L\int \int_{r\in  B_{\sqrt{\rho\nu}}} \left( \frac{R_mr^{2m+2}}{2\nu^m  \rho^m \sigma^m}\right) \frac{1}{\sqrt{4\pi \rho\nu \sigma}^D}\exp\left(-\frac{r^2}{4\rho\nu \sigma} \right)drd\sigma {\Big |},
\end{array}
\end{equation}
where the rapidly decreasing series of numbers $R_k$ is recursively defined by
\begin{equation}
R_1=\frac{1}{5},~R_{k+1}=\frac{1}{2}\frac{1}{6+2k-1}R_k,~\mbox{for all $k\geq 1$}.
\end{equation}
%
We observe that on a small time interval $[0,\Delta]$ we have a converging alternating series on the right sight of (\ref{series3}), because the members of the series behave like $O\left(\frac{1}{\sqrt{k}} \right)$ - as is easily inferred from Stirlings formula etc., i.e., the numbers $R_k$ converge fast enough such that the right side of the inequality
\begin{equation}\label{series4}
\begin{array}{ll}
{\Big |}\rho L\int \int_{\left\lbrace y|  y_i\geq 0\right\rbrace \cap B_{\sqrt{\rho\nu}}} \left( \frac{4y^2_i}{4\nu  \rho \sigma}\right) G^{\rho}_{\nu}(\sigma,y)dyd\sigma {\Big |}\\
\\
\leq \sum_{k=1}^{\infty}{\Big |}\rho L\int   \frac{R_k}{\sigma^k} \frac{1}{\sqrt{4\pi \sigma}^D}\exp\left(-\frac{1}{4\sigma} \right)d\sigma {\Big |}
\end{array}
\end{equation}
is finite, and the absolute value is decreasing (a detailed analysis than in this remark will be provided elsewhere). 
Hence for small $\sqrt{\rho\nu}$ we surely have the upper bound
\begin{equation}
\begin{array}{ll}
{\Big |}\rho L\int \int_{\left\lbrace y| y_i\geq 0\right\rbrace \cap B_{\sqrt{\rho\nu}}} \left( \frac{4y^2_i}{4\nu  \rho \sigma}\right) G^{\rho}_{\nu}(\sigma,y)dyd\sigma {\Big |}\\
\\
\leq {\Big |}\rho L\int   \frac{1}{4\sigma} \frac{1}{\sqrt{4\pi \sigma}^D}\exp\left(-\frac{1}{4\sigma} \right)d\sigma {\Big |}.
\end{array}
\end{equation}
For a small time step size $\Delta$ this upper bound becomes small as
\begin{equation}\label{parafreeest}
 {\Big |}\rho L\int_{0}^{\Delta}   \frac{1}{4} \frac{1}{\sqrt{4\pi}^D}\frac{1}{\sigma^{2.5}}\exp\left(-\frac{1}{4\sigma} \right)d\sigma {\Big |}
\end{equation}
becomes small for small $\Delta >0$.
Note that this upper bound becomes small as $\Delta$ becomes small independently of the parameters $\nu$ or $\rho$. Note that this additional scaling is obtained by a) the Lipschitz continuity of the nonlinear terms in the probabilistic local time solution representation, and b) by the combination of this Lipschitz continuity with the possibility to rewrite the equation in terms of convolutions with first order derivatives of the Gaussian. An analogous estimate as in (\ref{parafreeest}) holds oaslo on the interval $[\tau_{\rho 0},\tau_{\rho 0}+\Delta]$ (assuming that we have regular $H^m\cap C^m$ data for $m\geq 2$ at time $\tau_{\rho 0}$).
The estimate in (\ref{parafreeest}) implies that 
\begin{equation}\label{parafreeest2}
 {\Big |}\rho L\int_{0}^{\Delta}   \frac{1}{4} \frac{1}{\sqrt{4\pi}^D}\frac{1}{\sigma^{2.5}}\exp\left(-\frac{1}{4\sigma} \right)d\sigma {\Big |}\in o(\Delta)~ \mbox{for small $\Delta$},
\end{equation}
which implies in turn that a potential damping term introduced by a time transformation (cf. below)  becomes dominant if a regular norm exceeds a certain level.  
 For a local solution with data as in (\ref{vrhoinit}) you can check that $B^0_i, L^0,~B_i, L_i,~ B_{i,j}, L_{i,j}$ are Lipschitz for $1\leq i,j\leq D$. These considerations together with application of Young inequalities lead to the conclusion that forall $\tau_{\rho}\in \left[\tau_{\rho 0},\tau_{\rho 0}+\Delta_0 \right]$ we have
\begin{equation}
\max_{1\leq i\leq D}{\big |}B^0_i\ast G^{\rho}_{\nu,k}(\tau_{\rho},.){\big|}_{H^2\cap C^2}+{\big|}L^0\ast G^{\rho}_{\nu,k}(\tau_{\rho},.){\big |}_{H^2\cap C^2}\leq \rho c_DC^2_mC_KM_2
\end{equation}
or for all $1\leq i\leq D$
\begin{equation}
{\big |}B_i\ast G^{\rho}_{\nu}(\tau_{\rho},.){\big|}_{H^2\cap C^2}+{\big|}L_i\ast G^{\rho}_{\nu,k}(\tau_{\rho},.){\big |}_{H^2\cap C^2}\leq \rho c_DC^2_mC_KM_2,
\end{equation}
where $c_D$ is a finite constant which depends only on dimension $D$, and where $m=2$ may be chosen. This is an upper bound of growth which can be offset by an autocontrol scheme.
We may choose$\Delta >0$ small such that $M_2$ becomes small (smaller than $1$ especially) and then choose any $\rho$ with
\begin{equation}\label{rho}
0<\rho \leq \frac{1}{4c_DC^2_mC_K}
\end{equation}
in order to obtain a global scheme.
For the Navier Stokes equation, i.e., for fixed $\nu >0$ it is possible to analyze the damping term and compare it with upper bounds of the increment of the nonlinear terms over the time interval $[\tau_{\rho 0},\tau_{\rho 0}+\Delta]$.
We cannot do this in the case of the  Euler equation, of course. However the upper bound constructed above indicates that we can set up an autocontrol scheme which works without viscosity damping estimates. Here it is essential to have estimates which are independent of $\nu >0$, i.e.,  estimates which work in the case of the Navier Stokes equation  uniformly with respect to the viscosity $\nu >0$ (we then may invoke Lipschitz continuity of the convoluted Leray projection term and Burgers term, and apply some compactness argument in order to obtain a  viscosity limit as in the main text). So let us consider fixed $\nu >0$ in this remark. 

For $\tau_{\rho 0}\geq 0$ and a time interval $\left[\tau_{\rho 0},\tau_{\rho 0}+\Delta_0 \right]$ for some $\Delta_0 \in (0,1)$ we consider a comparison function $u^{\rho,\nu,\tau_0}_i,~1\leq i\leq D$ by a global time transformation with local time dilatation on a short time interval, i.e., we consider the transformation
\begin{equation}\label{utr1}
\begin{array}{ll}
(1+\tau_{\rho} )u^{\rho,\nu,\tau_{\rho 0}}_i(s,y)=v^{\rho,\nu}_i(\tau_{\rho},y),\\
\\
\mbox{where}~s=\frac{\tau_{ \rho}-\tau_{\rho 0}}{\sqrt{1-(\tau_{\rho}-\tau_{\rho 0})^2}},~\tau_{\rho}-\tau_{\rho 0}\in [0,\Delta_0],~1\leq i\leq D
\end{array}
\end{equation}
where $\Delta_0 \in [0,1)$. 
Using the abbreviation $\Delta \tau_{\rho}=\tau_{\rho}-\tau_{\rho 0}$, we have
\begin{equation}\label{utr2}
 v^{\rho,\nu}_{i,\tau_{\rho}}=
u^{\rho,\nu,\tau_{\rho 0}}_i
+(1+\tau_{\rho} )u^{\rho,\nu,\tau_{\rho 0}}_{i,s}\frac{ds}{d\tau_{\rho}},~\frac{ds}{d\tau_{\rho}}
=\frac{1}{\sqrt{1-\Delta \tau^2_{\rho}}^3}.
\end{equation}
We choose a time interval length $\Delta_0 \leq 0.5$ such that local time contraction, and, hence, a local time solution exists on the interval $[\tau_{\rho 0},\tau_{\rho 0}+\Delta_0]$. The equation for $u^{\rho,\nu,\tau_{\rho 0}}_i,~1\leq i\leq D$ becomes
\begin{equation}
\begin{array}{ll}\label{utr3}
u^{\rho,\nu,\tau_{\rho 0}}_{i,s}+\frac{\sqrt{1-\Delta \tau_{\rho}^2}^3}{1+\tau_{\rho}} \rho \nu \Delta u^{\rho,\nu,\tau_{\rho 0}}_i
+ \rho\frac{\sqrt{1-\Delta \tau^2_{\rho}}^3}{1+\tau_{\rho}}\sum_{j=1}^D \left( u^{\rho,\nu,\tau_{\rho 0}}_j\frac{\partial u^{\rho,\nu,\tau_{\rho 0}}_i}{\partial x_j}\right) \\
\\
-\rho \frac{\sqrt{1-\Delta \tau^2_{\rho}}^3}{1+\tau_{\rho}} \int_{{\mathbb R}^D} K_{D,i}(.-y)\sum_{j,m=1}^D\left( \frac{\partial u^{\rho,\nu,\tau_{\rho 0}}_m}{\partial x_j}\frac{\partial u^{\rho,\nu,\tau_{\rho 0}}_j}{\partial x_m}\right) (.,y)dy\\
\\
-\frac{\sqrt{1-\Delta \tau^2_{\rho}}^3}{1+\tau_{\rho}}u^{\rho,\nu,\tau_{\rho 0}}_i=0,\\
\\
~u^{\rho,\nu,\tau_{\rho 0}}(0,.)=v^{ \rho,\nu}(\tau_{\rho 0},.).
\end{array}
\end{equation}
This equation is defined on a dilated time interval $\left[0,\Delta_{d} \right] $, where the relation of $\Delta_d$ to $\Delta_0$ follows from the time transformation.
Now assume that an arbitrary large but fixed time horizon $T_{\rho}>0$ is given. Note that this corresponds to an original time horizon $T=\rho T_{\rho}$ with $\rho$ as in (\ref{rho}). Since $\rho$ depends only on the dimension and the data size in the $H^m\cap C^m$-norm essentially, a spatial norm of the initial data (and the spatial kernel constant $C_K$) global results with arbitrary horizon $T_{\rho}$ transfer to global results with respect to original time. 
Assume that for $m\geq 2$ and for some $0\leq \tau_{\rho 0}<T$ we have a global upper bound
\begin{equation}
\max_{1\leq i\leq D}{\Big |}v^{\rho,\nu}_i(t_0,.){\Big |}_{H^m\cap C^m}\leq C(1+\tau_0)<\infty.
\end{equation}
We may assume that $C\geq 3$. Then 
\begin{equation}
\max_{1\leq i\leq D}{\Big |}u^{\rho,\nu,\tau_{\rho 0}}_i(0,.){\Big |}_{H^m\cap C^m}=\frac{1}{1+\tau_{\rho 0}} \max_{1\leq i\leq D}{\Big |}v^{\rho,\nu}_i(\tau_{\rho 0},.){\Big |}_{H^m\cap C^m}\leq C,
\end{equation}
if $\rho$ is chosen as in (\ref{rho}).
Hence,
\begin{equation}
\max_{1\leq i\leq D}{\Big |}u^{\rho,\nu,\tau_{\rho 0}}_i(\Delta_d,.){\Big |}_{H^m\cap C^m}\leq C.
\end{equation} 
This implies that 
\begin{equation}
\begin{array}{ll}
\max_{1\leq i\leq D}{\Big |}v^{\rho,\nu}_i(\tau_{\rho 0}+\Delta_0,.){\Big |}_{H^m\cap C^m}\\ 
\\
\leq \max_{1\leq i\leq D}(1+\tau_{\rho 0}+\Delta_0){\Big |}u^{\rho,\nu,\tau_{\rho 0}}_i(\Delta_d,.){\Big |}_{H^m\cap C^m}
 \leq C(1+\tau_{\rho 0}+\Delta_0).
\end{array}
\end{equation} 
Iterating this procedure for all integers $k\geq 0$ with $k\Delta_0 \leq T_{\rho}$ we have 
\begin{equation}
\begin{array}{ll}
\max_{1\leq i\leq D}{\Big |}v^{\rho,\nu}_i(k\Delta_0,.){\Big |}_{H^m\cap C^m}\leq \max_{1\leq i\leq D}(1+k\Delta_0){\Big |}u^{\rho,\nu,\tau_{\rho 0}}_i(0,.){\Big |}_{H^m\cap C^m}\\
\\
\leq C(1+k\Delta_0).
\end{array}
\end{equation}
Then by a local time contraction we can interpolate and have for all $0\leq\tau_{\rho}\leq T_{\rho}$
\begin{equation}
\begin{array}{ll}
\max_{1\leq i\leq D}{\Big |}v^{\rho,\nu}_i(\tau_{\rho},.){\Big |}_{H^m\cap C^m}\leq \max_{1\leq i\leq D}(1+\tau_{\rho}){\Big |}u^{\rho,\nu,\tau_{\rho 0}}_i(0,.){\Big |}_{H^m\cap C^m} \\
\\
\leq 2C(1+\tau_{\rho}).
\end{array}
\end{equation}
\end{rem}

\begin{rem}
The argument below shows that singularities can be constructed even for data $\omega^f_i\in C^{\infty}_0,~1\leq i\leq D$, where $C^{\infty}_0$ the space of smooth functions vanishing at infinity. However, $\omega^f_i\in H^2\cap C^2,~1\leq i\leq D$ is essential, and we consider further regularity results and long-time kinks and singularities elsewhere. Here we are interested mainly in the connection to singular solution or kinks of the Navier Stokes equation with regular time dependent force terms.  
\end{rem}
%
Theorem \ref{main1} can be extended in the sense that there are solution branches with weak singularities of any integer order $k$ in the sense that a solution function is only in $C^{k-1}\setminus C^{k}$. In order to have a succinct statement we introduce the concept of a 'spatial kink of order $k$' of a classical solution.

\begin{defi}
We say that a classical solution branch $\omega_i,~1\leq i\leq D$ of an incompressible Euler Cauchy problem with data $\omega^{f,-}_{i},~1\leq i\leq D$ with  $\omega^f_i\in C^{m}$ for $m\geq k$ has spatial a kink of order $k\geq 1$, if there is a space-time point $(\tau,x)$ with $\tau>0$ and $x\in {\mathbb R}^D$ such that
\begin{equation}
\omega_i\in C^{k-1}\setminus C^{k}~\mbox{at }~(\tau,x)
\end{equation}
for some integer $k\geq 1$.
\end{defi}

\begin{cor}\label{maincor}
Let $D=3$. For any $k\geq 2$ and $s\geq 0$ there exist data $\omega^f_i\in H^m\cap C^{m},~1\leq i\leq D$ with $m\geq k+2$ and a vorticity solution $\omega_i,~1\leq i\leq D$ of the three dimensional incompressible Euler equation Cauchy problem with data $\omega^f_i,~1\leq i\leq D$ such that after some finite time the solution has a kink of order $k$.
\end{cor}

\begin{rem}
We treat the latter result as a Corollary because the proof method is the same, and only the choice of data for the time reversed Euler-type Cauchy problem has to be adapted. 
\end{rem}

Next we draw consequences for the Navier Stokes equation. A classical solution 
\begin{equation}
\omega_i,~1\leq i\leq D,~\omega_i\in C^{1}\left(\left[0,T\right), H^2\cap C^2 \right) 
\end{equation}
of the incompressible Euler equation (on the interval $[0,T)$ without the point $T$) and  with a blow-up of vorticity at time $T$ satisfies the vorticity form of the Navier Stokes equation
 \begin{equation}\label{vorticitynav}
\frac{\partial \omega}{\partial \tau}-\nu \Delta \omega+v\cdot \nabla \omega=\frac{1}{2}\left(\nabla v+\nabla v^T\right)\omega +F, 
\end{equation}
with force term $F=\left(F_1,F_2,F_3 \right)^T$ (on the same time interval $[0,T)$) if
\begin{equation}\label{Fterm}
\mbox{for all}~t~F_i(t,.)=-\nu \Delta \omega_i(t,.)\in L^2\cap C,~1\leq i\leq 3.
\end{equation}
The analysis below shows that $F_i$ is also $L^2$ with respect to time on the time interval $[0,T]$, where $T>0$ is the time where the vorticity of the Euler equation blows up.  
Similarly, a regular classical solution on the time interval $[0,T)$
\begin{equation}
\omega_i,~1\leq i\leq D,~\omega_i(\tau,.)\in H^{2+k}\cap C^{2+k}
\end{equation}
of the incompressible Euler equation with a kink of order $k$ at time $T>0$ such that
\begin{equation}
\omega(T,.)\in C^{k-1}\setminus C^k,
\end{equation}
satisfies the vorticity form of the Navier Stokes equation (\ref{vorticitynav}),
where
\begin{equation}\label{Fterm}
\mbox{ for all $t\in [0,T)$}~~F_i=-\nu \Delta \omega_i(\tau,.)\in H^k\cap C^k,~1\leq i\leq 3.
\end{equation}
We conclude
\begin{thm}
For the Cauchy problem for the Navier Stokes equation (\ref{vorticitynav}) for some time $T>0$ time dependent force terms with
\begin{equation}
F_i(\tau,.)\in H^k\cap C^k,~k\geq 0,~\mbox{for all }~\tau\in[0,T)
\end{equation}
and data  $\omega_i(0,.)\in H^{2+k}\cap C^{2+k},~1\leq i\leq 3,~k\geq 0$ can be chosen such that a regular classical solution of the Navier Stokes equation on the time interval $[0,T)$ has a blow-up (case $k=0$) or can be extended beyond the time interval $[0,T]$ and then has a kink of order $k\geq 1$ at time $T$.
\end{thm}


\section{Proof of Theorem \ref{main1} and Corollary \ref{maincor}}
For the Gaussian fundamental solution $G_{\nu}$ of the equation 
\begin{equation}
p_{,t}-\nu \Delta p=0
\end{equation}
we consider on some time interval $[0,T]$ with time horizon $T>0$ the iteration scheme $v^{\nu,-,k}_i,~1\leq i\leq D,~k\geq 0$, where for $k\geq 1$
\begin{equation}\label{Navlerayscheme}
\begin{array}{ll}
 v^{\nu,-,k}_i=v^{f,-}_i\ast_{sp}G_{\nu}
+\sum_{j=1}^D \left( v^{\nu,-,k-1}_j\frac{\partial v^{\nu,-,k-1}_i}{\partial x_j}\right) \ast G_{\nu}\\
\\-\left( \sum_{j,m=1}^D\int_{{\mathbb R}^D}\left( \frac{\partial}{\partial x_i}K_D(.-y)\right) \left( \frac{\partial v^{\nu,-,k-1}_m}{\partial x_j}\frac{\partial v^{\nu,-,k-1}_j}{\partial x_m}\right) (.,y)dy\right) \ast G_{\nu},\\
\\
\mbox{and where this scheme is initialized by}\\
\\
v^{\nu,-,0}_i:=v^{f,-}_i(.)\ast_{sp}G_{\nu}.
\end{array}
\end{equation}
Here the velocity data $v^{f,-}_i,~1\leq i\leq D$ are determined by the vorticity data $\omega^{f,-}_i,~1\leq i\leq D$ via the Biot-Savart law. The upper script $^-$ reminds us that we consider a time-reversed problem auch that the data chosen here are not confused with the initial data $\omega^{f}_i,~1\leq i\leq D$ in the statement of our main theorem (which are the final data at some time $T>0$ of the time reversed problem).   In this paper we choose a refined version  of the scheme in (\ref{Navlerayscheme}), where we  construct a singularity for at least one velocity component.  We single out one velocity component in $H^1\setminus C^1$, and initialize the scheme in (\ref{Navlerayscheme}) with
\begin{equation}
v^{\nu,-,0}_{i_0}:=v^{f,-}_{i_0}(.)\ast_{sp}G_{\nu}
\end{equation}
as before, while for $j\neq i_0$ we choose
\begin{equation}
v^{\nu,-,0}_{j}:=v^{f,-}_{j}(.),
\end{equation}
where we observe below that the data $v^f_{j}(.),~j\neq i_0$ can be chosen to have full regularity (cf. item i) below). We shall consider possible constructions of data with stronger and weaker singularities below, where we also consider some variations of argument depending on the choice of the data.
 
Note that we  start the iteration scheme with smoothed data $v^{f,-}_{i_0}(.)\ast_{sp}G_{\nu}$ which appear in the first approximating increment $\delta v^{\nu,1,-}_i=v^{\nu,-,1}_i-v^{f,-}_i\ast_{sp}G_{\nu},~1\leq i\leq D$. This makes it possible to use estimates of $v^{f,-}_{i_0}\ast_{sp}G_{\nu,j},~1\leq i,j\leq D$ in the estimation of the functional increments of the iterated scheme. In a former version of this paper we started the iteration scheme with the data $v^{f,-}_i,~ 1\leq i\leq D$ themselves, and this has the effect that the smoothing in the increment is felt only after the second iteration. Smoothing coefficients is often useful (sic Leray). Here we use smoothed approximations in the nonlinear terms of local iteration schemes.

 We note that the symbol $\ast_{sp}$ denotes convolution with respect to the spatial variables, and $\ast$ denotes convolution with respect to the time variable and the spatial variables. It is an advantage to start the scheme with the spatially convoluted data $v^{f,-}_{i}(.)\ast_{sp}G_{\nu}$ at least for one velocity component, since we choose the data $v^{f,-}_i,~1\leq i\leq D$ to be only locally Lipschitz for one velocity component $v_{i_0}$. The approximation of first order derivatives $v^{f,-}_i\ast_{sp}G_{\nu,j},~1\leq i,j\leq D$ in the nonlinear functional increments $\delta v^{\nu,1,-}_i$ is an advantage for the estimates. Note that all velocity component data are chosen such that they are smooth in the complement of the origin and have strong spatial decay at spatial infinity. 
\begin{rem} 
Note that the data we choose below are only locally Lipshitz and globally H\"older continuous. If we apply the argument to globally Lipschitz continuous data ( e.g. $\alpha_0=\beta_0$ in the choice of data in item i) below), then we obtain weaker singularities of bounded oscillatory type.
\end{rem}
  
All data for all velocity components are in $H^1$ such that the first derivatives in the approximation of the velocity have the representation
\begin{equation}
v^{\nu,-,0}_{i,j}:=v^{f,-}_i(.)\ast_{sp}G_{\nu,j}=v^{f,-}_{i,j}(.)\ast_{sp}G_{\nu}.
\end{equation} 
The data $v^{\nu,-,0}_{i},~ 1\leq i\leq D$ themselves are H\"{o}lder continuous with exponent $\beta_0$ close to $1$ (even Lipschitz at the origin) for one velocity component and, depending on the variation of arguments considered below, they are also H\"older continuous, or  Lipschitz continuous, or even of full regulrity for the other velocity components.  
We shall use the fact that Lipschitz continuous data convoluted with the first order derivatives of the Gaussian have convenient upper bounds, which are even $\nu$-independent. However, this is not decisive since we have this effect for $v^{\nu,-,k}_i$ for $k\geq 2$ anyway if we start the scheme with $v^{f,-}_i(.)$ instead of $v^{f,-}_i(.)\ast_{sp}G_{\nu}$ at $k=0$.
We choose an index $i_0\in \left\lbrace 1,\cdots,D\right\rbrace$ and data $v^f_{i_0}$ and complement data $v^f_{j},~j\in \left\lbrace 1,\cdots,D\right\rbrace\setminus \{i_0\}$ such that corresponding vorticity component $\omega^{f,-}_{i_0}$ are singular at one point. Furthemore the complement data $v^{f,-}_{j},~j\in \left\lbrace 1,\cdots,D\right\rbrace\setminus \{i_0\}$ are constructed with the iteration scheme $v^{\nu,-,k}_i,~1\leq i\leq D,~k\geq 1$ such that for all $k$
\begin{equation}\label{incompk}
\mbox{ for all $t\in (0,T]$}~~\sum_{i=1}^Dv_{i,i}=\lim_{k\uparrow \infty}\sum_{i=1}^Dv^{\nu,-,k}_{i,i}(t,.)=0.
\end{equation}
The basic idea to realise the additional constraint in (\ref{incompk}) is as follows. Local contraction leads to a local representation
\begin{equation}\label{Navlerayschemefix}
\begin{array}{ll}
 v^{\nu,-}_i=v^{f,-}_i\ast_{sp}G_{\nu}
+\sum_{j=1}^D \left( v^{\nu,-}_j\frac{\partial v^{\nu,-}_i}{\partial x_j}\right) \ast G_{\nu}\\
\\-\left( \sum_{j,m=1}^D\int_{{\mathbb R}^D}\left( \frac{\partial}{\partial x_i}K_D(.-y)\right) \left( \frac{\partial v^{\nu,-}_m}{\partial x_j}\frac{\partial v^{\nu,-}_j}{\partial x_m}\right) (.,y)dy\right) \ast G_{\nu}\\
\\
=:v^{f,-}_i\ast_{sp}G_{\nu}+\delta v^{\nu,-}_i.
\end{array}
\end{equation}
Applying the divergence operator we get for any $\nu >0$
\begin{equation}\label{Navlerayschemefix}
\begin{array}{ll}
\sum_{i=1}^D v^{\nu,-}_{i,i}=\sum_{i=1}^Dv^{f,-}_i\ast_{sp}G_{\nu ,i}
+\sum_{i,j=1}^D \left( v^{\nu,-}_{j,i} v^{\nu,-}_{i,j}\right) \ast G_{\nu}\\
\\-\sum_{i=1}^D\left( \sum_{j,m=1}^D\int_{{\mathbb R}^D}\left( \frac{\partial^2}{\partial x_i^2}K_D(.-y)\right) \left( \frac{\partial v^{\nu,-}_m}{\partial x_j}\frac{\partial v^{\nu,-}_j}{\partial x_m}\right) (.,y)dy\right) \ast G_{\nu}.
\end{array}
\end{equation}
In a viscosity limit $\nu_k\downarrow 0$ the pressure condition
\begin{equation}
\Delta p=\lim_{k\uparrow \infty}\sum_{i,j=1}^D \left( v^{\nu_k,-}_{j,i} v^{\nu_k,-}_{i,j}\right)
\end{equation}
is satisfied, where 
\begin{equation}
-\sum_{i=1}^D \sum_{j,m=1}^D\int_{{\mathbb R}^D}\left( \frac{\partial^2}{\partial x_i^2}K_D(.-y)\right) 
\sum_{j,m=1}^D\left( \frac{\partial v^{\nu,-}_m}{\partial x_j}\frac{\partial v^{\nu,-}_j}{\partial x_i}\right) =-\Delta p.
\end{equation}
Hence, if the increment $\lim_{k\uparrow \infty}\delta v^{-,\nu_k}_i,~1\leq i\leq D$ is regular, then the incompressibility condition (on a time interval $(0,T]$) for a viscosity limit 
\begin{equation}
\lim_{k\uparrow \infty}\sum_{i=1}^D v^{\nu_k,-}_{i,i}=0 
\end{equation}
becomes equivalent to
\begin{equation}\label{vvisc}
\lim_{k\uparrow \infty}\sum_{i=1}^D v^{\nu_k,-}_{i,i}=\lim_{k\uparrow \infty}\sum_{i=1}^Dv^{f,-}_i\ast_{sp}G_{\nu_k ,i}=0
\end{equation}
on a time interval $(0,T]$. 
Now note that for our choice of data $v^f_{i_0}$ below at item i) we have for given $\nu >0$ and $t>0$
\begin{equation}
(v^{f,-}_{i_0}\ast_{sp}G_{\nu,i_0})(t,0)=\int \phi_0(r)\frac{-(0-y_{i_0})}{4\nu t}G_{\nu}(t,0-y)dy=0,
\end{equation}
since the data $\phi_0$ (cf. below) are symmetric (here $r=\sqrt{\sum_{i=1}^Dy_i^2}$). We shall use data below which are locally regular in the complement of the origin. Next we observe for any global bounded continuous and locally regular function $f$   we have a degeneracy in the viscosity limit as with $y_i=\nu t z_i,~1\leq i\leq D$
\begin{equation}\label{degen}
\begin{array}{ll}
{\Big |}\int f(x-y)\frac{-(y_{i})}{4\nu t}G_{\nu,i}(t,y)dy{\Big |}={\Big |}\int f(x-y)\frac{-(y_{i})}{4\nu t}\frac{1}{\sqrt{4\pi\nu t}^D}\exp\left(-\frac{|y|^2}{4\nu t}\right) dy{\Big |}\\
\\
\leq {\Big |}\int f(x-\nu tz)z_i\frac{1}{\sqrt{4\pi\nu t}^D}\exp\left(-|z|^2\nu t\right) (\nu t)^Ddz{\Big |}\downarrow 0~\mbox{as $\nu\downarrow 0$.}
\end{array}
\end{equation}
Here local regularity around $x$ and a local multivariate Taylor formula can be used. More precisely, we may use for $m=1$
\begin{lem}
Let $x\in U\subset {\mathbb R}^D$  be an open set. 
For $1\leq m\geq l\geq 0$ and $f\in C^m(U)$
 we have for all $x,h\in {\mathbb R}^D$ with $\{x +sh|s\in [0,1]\}\subset U$
\begin{equation}
\begin{array}{ll}
f(x+h)=f(x)+\sum_{0<|\alpha| <m}\frac{(\partial^{\alpha}f)(x)}{\alpha !}\\
\\
+m\sum_{|\alpha|=m}\frac{h^{\alpha}}{\alpha !}\int_0^1(1-\theta)^{m-1}\left(\partial^{\alpha}f \right)(x+\theta h)d\theta ,
\end{array}
\end{equation}
where $\partial^{\alpha}f$ denotes the multivariate partial derivative of $f$ of order $\alpha$ with respect to the multiindex $\alpha=(\alpha_1,\cdots,\alpha_D)$.
\end{lem}
In (\ref{degen}) the right term of the inequality can be estimated by two summands. The first summand of order zero in $\nu tz$ is
\begin{equation}
{\Big |}\int f(x)z_i\frac{1}{\sqrt{4\pi\nu t}^D}\exp\left(-|z|^2\nu t\right) (\nu t)^Ddz{\Big |}\downarrow 0~ \mbox{as}~ \nu \downarrow 0.
\end{equation}
Furthermore, if for $0\neq x$ $f$ is locally regular, then there is a neighborhood $U$ of $x$ where the first order partial derivatives are bounded. If $z$ is in any compact (large) neighborhood then $\{x+\theta \nu tz,~ \theta \in[0,1]\}$ is in $U$ and the additional fcato $\nu tz $ of the Taylor formula remainder term ensures that the integral over any compact set converges to zero as $\nu$ converges to zero. It is then observed straightforwoardly that the integral over the complement of that large compact set also converges to $0$ as $\nu$ converges to $0$.
  
However in the local solution increments $\delta v^{\nu,-,k}_i,~1\leq i\leq D$ the first order derivatives of the data
\begin{equation}
v^{\nu,-,0}_{i_0,k}=v^{f,-}_{i_0,k}\ast_{sp}G_{\nu},~v^{\nu,-,0}_{j,k}=v^{f,-}_{j,k},~j\neq i_0
\end{equation}
for all $1\leq k\leq D$ contribute to the local solution, and this also holds in the viscosity limit. 

 We note that especially for given $\nu >0$  $t\in (0,T]$  $x\rightarrow (v^{f,-}_{i_0}\ast_{sp}G_{\nu,i_0})(t,x)$ is continuous and bounded by a constant which is independent of $\nu$ and $t$. 
Hence our analysis of the functional solution increment below shows that we may choose functions $v^{f,-}_{j}$ for $j\neq i_0$ which are $C^1$ and have bounded first order derivatives such that the incompressibility condition
\begin{equation}
v^{f,-}_{i_0}\ast_{sp}G_{\nu,i_0}+\sum_{j\neq i_0}v^{f,-}_{j}\ast_{sp}G_{\nu,j}=0
\end{equation}
is satisfied on $(0,T]$ and such that the incompressibility relation  $\sum_{i=1}^Dv^{-}_{i,i}=\lim_{k\uparrow \infty}\sum_{i=1}^D v^{\nu_k,-}_{i,i}$ holds pointwise on $(0,T]$ in the viscosity limit. Well for one variation of argument and our choice below of the data $v^{f,-}_{i_0}$ we can even choose data $v^{f,-}_j$ in Schwartz space for $j\ne i_0$. We can make the general statement: the data $v^{f,-}_i,~i\neq i_0$ can be chosen in the same regularity class (ar least) as the solution increment $\delta v^{-}_i=v^{-}_i-v^{f,-}_i,~1\leq i\leq D$ of the local time solution increment of the Euler equation.

Concerning the regularity of the local solution increments consider their first approximations in the iteration scheme. We have for $\nu >0$
\begin{equation}\label{Navlerayschemefix2}
\begin{array}{ll}
\delta v^{\nu,-,1}_i =v^{\nu,-,1}_i-v^{f,-}_i\ast_{sp}G_{\nu}\\
\\
=\sum_{j=1}^D \left( \left( v^{f,-}_j\ast_{sp}G_{\nu}\right) \left( v^{f,-}_i\ast_{sp}G_{\nu,j}\right) \right) \ast G_{\nu}\\
\\-{\Big (} \sum_{j,m=1}^D\int_{{\mathbb R}^D}\left( \frac{\partial}{\partial x_i}K_D(.-y)\right)\times \\
\\
\times \left( \left( v^{f,-}_m\ast_{sp}G_{\nu,j}\right)  \left( v^{f,-}_j\ast_{sp}G_{\nu,m}\right) \right) (.,y)dy{\Big )} \ast G_{\nu}\\
\\
=:v^{f,-}_i\ast_{sp}G_{\nu}+\delta v^{\nu,-,1}_i.
\end{array}
\end{equation}
Only one component of the velocity data, i.e., $v^{f,-}_{i_0}$ is chosen such that the vorticity has a singularity; the other are chosen to be regular. The regularity of the other data $v^{f,-}_j,~j\neq i_0$ can be exploited in order to represent second order spatial derivatives  $v^{\nu,-}_{i_0,i_0,i_0}$ by
\begin{equation}
v^{\nu,-}_{i_0,i_0,i_0}=-\sum_{j\neq i_0}v^{\nu,-}_{j,j,i_0}.
\end{equation}
For any $\nu >0$ we get the local representation 
\begin{equation}\label{Navlerayschemefix3}
\begin{array}{ll}
\delta v^{\nu,-,1}_{i_0,i_0,i_0} =v^{\nu,-,1}_{i_0,i_0,i_0}-v^{f,-}_{i_0,i_0}\ast_{sp}G_{\nu,i_0}\\
\\
=\sum_{j\neq i_0}\left( \left( v^{f,-}_j\ast_{sp}G_{\nu}\right) \left( v^{f,-}_{i_0}\ast_{sp}G_{\nu,j}\right) \right)_{,i_0} \ast G_{\nu,i_0}\\
\\
+\left( \left( v^{f,-}_{i_0}\ast_{sp}G_{\nu}\right) \left(- \sum_{j\neq i_0}v^{f,-}_{j}\ast G_{\nu,j}\right) \right)_{,i_0} \ast G_{\nu,i_0}\\
\\- \sum_{j,m\neq i_0}^D\int_{{\mathbb R}^D}K_{D,i_0}(.-y)\times\\
\\
\times \left( \left( v^{f,-}_m\ast_{sp}G_{\nu,j}\right)  \left( v^{f,-}_j\ast_{sp}G_{\nu,m}\right) \right)_{,i_0}(.,y)dy \ast G_{\nu,i_0}\\
\\
-2\sum_{j=1}^D\int_{{\mathbb R}^D}K_{D,i_0}(.-y) 
 \left( \left( v^{f,-}_{i_0}\ast_{sp}G_{\nu,j}\right) 
\left( v^{f,-}_j\ast_{sp}G_{\nu,i_0}\right) \right)
 (.,y)dy\ast G_{\nu,i_0},
\end{array}
\end{equation}
where in the last line $\left( v^{f,-}_{i_0}\ast_{sp}G_{\nu,j}\right)$ can be substituted by $\sum_{j\neq i_0}v^{f,-}_{j}\ast G_{\nu,j}$ as in the Burgers term above.
Hence, as $v^{f,-}_j,~j\neq i_0$ are regular data, this can clearly be exploited in order to obtain regularity of second derivatives of the velocity component $v_{i_0}$. For our choice of data in item i) below we even have a further reduction (concerning regularity) for the terms $v^{f,-}_{i_0}\ast_{sp}G_{\nu,j}$.
Furthermore, we choose the data $v^{f,-}_i, 1\leq i\leq D$ such that Burgers term data and the Leray projection data in the iteration scheme are Lipschitz. Especially we choose the data such that the Burgers terms
\begin{equation}
\left( \left( v^{f,-}_j\ast_{sp}G_{\nu}\right) \left( v^{f,-}_i\ast_{sp}G_{\nu,j}\right) \right) 
\end{equation}
and the Leray projection terms
\begin{equation}
\int_{{\mathbb R}^D}\left( \frac{\partial}{\partial x_i}K_D(.-y)\right) 
\sum_{j,m=1}^D \left( \left( v^{f,-}_m\ast_{sp}G_{\nu,j}\right) 
 \left( v^{f,-}_j\ast_{sp}G_{\nu,m}\right)\right) (.,y)dy 
\end{equation}
are Lipschitz. Two remarks are in order here: a) we explain the role of Lipschitz continuity of the data, and b) we give the reason, why the construction works in dimension $D=3$ and not $D=2$ or $D=1$. 
Ad a), we shall choose the data $v^{f}_{i_0}$ to be locally Lipschitz continuous such that the symmetry of the first order spatial derivative of the Gaussian (cf. (\ref{symm}) and the related estimate below) such that
\begin{equation}
v^{\nu,-}_{i_0,i_0}(t,.)=v^{f,-}_i\ast_{sp}G_{\nu,i_0}(t,.)+\delta v^{\nu,-}_{i_0,i_0}(t,.)
\end{equation}
has an appropriate upper bound for $t\in (0,T]$, where regularity observations concerning the increment $\delta v^{\nu,-}_{i_0,i_0}(t,.)$ come into play. 
In the latter context we define for $k\geq 1$
\begin{equation}\label{vinitk}
\delta v^{\mbox{init},\nu,-,k}_i=v^{\nu,-,k}_i-v^{f,-}_i\ast_{sp}G_{\nu},
\end{equation}
and
\begin{equation}
\delta v^{\nu,-,k}_i=v^{\nu,-,k}_i-v^{\nu,-,k-1}_i.
\end{equation}
We show that for some time horizon $T>0$ the local time solution function $v^{\nu,-}_i,~1\leq i\leq D$ of the Navier Stokes type extension of the time-reversed Euler equation has a representation of the form
\begin{equation}\label{v-rep}
v^{\nu,-}_i=v^{f,-}_i\ast_{sp}G_{\nu}+\delta v^{\mbox{init},\nu,-,2}_i+\sum_{k=3}^{\infty}\delta v^{\nu,-,k}_i,
\end{equation}
which is spatially in $H^3\cap C^3$ for evaluations at all time $t\in (0,T]$.
\begin{rem}
Note the choice $k=2$ in (\ref{v-rep}) with respect to the initial increment ( cf. (\ref{vinitk})). However, as we start with the iteration scheme with the data $v^{\nu,-,0}_i:=v^{f,-}_i(.)\ast_{sp}G_{\nu}$ this is not necessary and the alternative representation
\begin{equation}\label{v-rep1}
v^{\nu,-}_i=v^{f,-}_i\ast_{sp}G_{\nu}+\delta v^{\mbox{init},\nu,-,1}_i+\sum_{k=3}^{\infty}\delta v^{\nu,-,k}_i,
\end{equation}
is also sufficient for this improved scheme. If we start with data $v^{f,-}_i(.)$ at $k=0$ (as we did in a former version of this paper)  instead, then the
 choice $k=1$ is not sufficient as the first initial increment in  (\ref{v-rep1}) is then  of lower regularity -at least for some of the possible choices of data considered below.
\end{rem}

Next ad b) we note that there is  a dependence on dimension in the construction of the viscosity limit. 
For one variation of argument where we use the degeneracy of of the first order derivative data term $v^{f,-}_{i_0}\ast_{sp}G_{\nu_k,i_0}$ in the viscosity limit for $\nu_k\downarrow 0$ this and $H^1\cap C^1$-regularity or even full reesgultity of the data function $v^{f,-}_{j}$ we need $D\geq 3$ obviously. 
A local solution of the time-reversed Euler equation is constructed via a viscosity limit of a subsequence $v^{\nu_k,-}_i,~ 1\leq i\leq D$ with $\nu_k\downarrow 0$ (a compactness argument). In this compactness argument we need a uniform upper bound
\begin{equation}
\max_{1\leq i,j\leq D}\sup_{0\leq t\leq T}{\big |}v^{f,-}_i\ast_{sp}G_{\nu_k,j}(t,.){\big | }\leq C
\end{equation}
for some finite constant $C>0$ which is independent of $\nu$ (or $\nu_k,~ k\geq 1$ at least).
Note that the classical Gaussian estimate
\begin{equation}\label{gaussestess}
\begin{array}{ll}
G_{\nu}(t,z)=\sqrt{\pi}^{-D/2}\left( \frac{|z|^2}{4\nu t}\right)^{D/2-\delta}\frac{1}{|z|^{D-2\delta}(4\nu t)^{\delta}}\exp\left( \frac{-|z|^2}{4\nu t}\right)\\
\\
\leq \frac{C}{|y|^{D-2\delta}(4\nu t)^{\delta}},
\end{array}
\end{equation}
where
\begin{equation}\label{constant}
C=\sup_{|w|>0}|w|^{D-2\delta}\exp(-|w|^2),
\end{equation}
imposes no serious restrictions regarding dimension. For first order spatial derivatives of the Gaussian we use Lipschitz continuous upper bounds of convoluted data terms, where we use spatial antisymmetry of the Gaussian /cf. below). For some variation of argument it can be useful to have $\delta >0.5$, such that the independence of the finite constant in (\ref{constant}) of $\nu$ is only ensured of $D>1$. However, we are interested in $D\geq 3$ only anyway. The $2$-dimensional flows with scalar vorticity lead to the constraint
\begin{equation}
\omega\cdot \nabla v=0~(\mbox{dim} D=2),
\end{equation}
(cf. \cite{MB}), which means that the Leray projection term cancels. Hence under this specific circumstances an inviscid Burgers equations results, which may also have singular solutions (the Burgers equation with in spatial dimension definitely has a shock wave). The vorticity form of the incompressible Euler- and Navier Stokes equation is clearly a strong constraint in case of diemsnion $D\leq 2$, such that it is no suprise that in case $D=2$ singularities can only be proved in case of more comlex problems with boundaries(cf. \cite{Ki}).

However, this is not our main interest, and we restrict our discussion to the case $D=3$ or $D\geq 3$ (as far as velocity is concerned). 
Let us pause for a moment here and show that the Leray projection operator and  the Burgers operator applied to data functions in $H^2\cap C^2$ are Lipschitz continuous functions which have Lipschitz continuous first order derivatives. We consider the Leray projection term (a simplified similar argument works for the Burgers term).
For $g=(g_1,\cdots,g_D)^T$ with $\max_{1\leq i\leq D}{\big |}g_i{\big |}_{ H^2\cap C^2}\leq C_2$ consider the function
\begin{equation}\label{Navlerayschemefix2}
\begin{array}{ll}
L_{gi}(.):=\sum_{j,m=1}^D\int_{{\mathbb R}^D}\left( \frac{\partial}{\partial x_i}K_D(.-y)\right) \left( \frac{\partial g_m}{\partial x_j}\frac{\partial g_j}{\partial x_m}\right) (y)dy.
\end{array}
\end{equation}
Here it is sufficient to consider functions $g_i,~ 1\leq i\leq D$  without time dependence, as the estimate transfer straightforwardly  to the local solution functions. 
For given $x,x'\in {\mathbb R}^D$ we verify Lipschitz continuity (thereby verifying Lipschitz continuity on compact domains). We have for $x,x'\in {\mathbb R}^D$
\begin{equation}
\begin{array}{ll}
{\big |}L_{gi}(x)-L_{gi}(x'){\big |}\leq  \\
\\
 {\big |}\sum_{j,m=1}^D\int_{{\mathbb R}^D}\left( K_{D,i}(x'-y)- K_{D,i}(x-y)\right) \left( \frac{\partial g_m}{\partial x_j}\frac{\partial g_j}{\partial x_m}\right) (y)dy{\big |} \\
 \\
 \leq  {\Big |}\sum_{j,m=1}^D\int_{{\mathbb R}^D}{\big |} K_{D,i}(x'-y)- K_{D,i}(x-y){\big |} {\Big |} C_2\frac{\partial g_j}{\partial x_m}{\Big |} (y)dy{\Big |}\\
 \\
 \leq  C_2{\Big |}\int_{{\mathbb R}^D}{\Big |}  K_{D,i}(y) {\Big |} 
 \sum_{j,m=1}^D {\Big |}  \frac{\partial g_j}{\partial x_m} (x'-y)-\frac{\partial g_j}{\partial x_m} (x-y){\Big |}  dy{\Big |}\\
 \\
 \leq  C_3{\Big |}\int_{{\mathbb R}^D}{\Big |}  K_{D,i}(y){\Big |}  
 \sum_{j,m,k=1}^D {\Big |} \int_0^1 g_{j,m,k} (x'-y+\theta(x-x'))d\theta|x'-x|{\Big |} dy{\Big |}\\
 \\
 \leq C_3C_Kl_{g1}|x-x'|,
\end{array}
\end{equation}
where we use the general Young inequality and the equivalence of finite dimensional norms, i.e. the equivalence of the Euclidean norm and the maximum norm of ${\mathbb R}^D$ (this leads to an additional constant which we absorbed into the constant $C_3$). Furthermore, where $l_{g1}$ is a maximum of Lipschitz constants of the first order spatial derivatives of the functions $g_j,~ 1\leq j\leq D$ and $C_K$ is an integral constant related to the Laplacian kernel. Note that all these constants are independent of $\nu$. It is clear that Lipschitz continuity on compact domains is sufficient for our purposes, but global Lipschitz continuity follows as well. Next for the first order spatial derivatives of the Leray projection term consider for $1\leq k\leq D$ the functions
\begin{equation}\label{Navlerayschemefix22}
\begin{array}{ll}
L_{gi,k}(.):=2\sum_{j,m=1}^D\int_{{\mathbb R}^D}\left( \frac{\partial}{\partial x_i}K_D(.-y)\right) \left( \frac{\partial g_m}{\partial x_j}\frac{\partial^2 g_j}{\partial x_m\partial x_k}\right) (y)dy.
\end{array}
\end{equation}
Concerning Lipschitz continuity of this function on compact domains we can argue as above. For given $x,x'\in {\mathbb R}^D$ we verify Lipschitz continuity (on compact domains). We have
\begin{equation}
\begin{array}{ll}
{\big |}L_{gi,k}(x)-L_{gi,k}(x'){\big |}\leq  \\
\\
 2{\big |}\sum_{j,m=1}^D\int_{{\mathbb R}^D}\left( K_{D,i}(x'-y)- K_{D,i}(x-y)\right) \left( \frac{\partial g_m}{\partial x_j}\frac{\partial^2 g_j}{\partial x_m\partial x_k}\right) (y)dy{\big |} \\
 \\
 \leq  2{\big |}\sum_{j,m=1}^D\int_{{\mathbb R}^D}{\Big |} K_{D,i}(x'-y)- K_{D,i}(x-y){\Big |}  {\Big |}  C_2\frac{\partial g_m}{\partial x_j}{\Big |}  (y)dy{\big |}\\
 \\
 \leq  2C_2{\big |}\sum_{j=1}^D\int_{{\mathbb R}^D}{\Big |}  K_{D,i}(y) {\Big |} 
  {\Big |}  \frac{\partial g_j}{\partial x_m} (x'-y)-\frac{\partial g_j}{\partial x_m} (x-y){\Big |}  dy{\big |}\\
 \\
 \leq 2C_3C_Kl_{g1}|x-x'|.
\end{array}
\end{equation}
Next note that the data $v^{f,-}_i,~1\leq i\leq D$ are $\nu$-independent H\"{o}lder continuous data such that we have
\begin{equation}\label{init0}
\lim_{\nu\downarrow 0}v^{f,-}_i\ast_{sp}G_{\nu}=v^{f,-}_i,~1\leq i\leq D.
\end{equation}
Moreover for Lipschitz continuous data $v^{f,-}_i$ we have even a pointwise upper bound for first order derivatives $v^{f,-}_i\ast_{sp}G_{\nu,j}$ which can be used in the iteration scheme. 
The convergence in (\ref{init0}) is in $H^2$ spatially (for the data chosen below). Hence we have H\"{o}lder continuous and pointwise convergence in case of dimension $D=3$. However, our main interest concerns the regularity of the increment
\begin{equation}
\delta v^{-}_i=\lim_{\nu\downarrow 0}v^{\nu,-}_i-v^{f,-}_i.
\end{equation}
For this purpose we need $\nu$-independent estimates of convolutions 
\begin{equation}
g_{\nu}\ast G_{\nu},~g_{\nu}\ast G_{\nu,i}
\end{equation}
with data
\begin{equation}
g_{\nu}\in S:=\left\lbrace ~\delta v^{\mbox{init},\nu,-,2}_i,~\delta v^{\nu,-,k}_i,~k\geq 3\right\rbrace .
\end{equation}
Note that $G_{\frac{\nu}{2}}$ is a probability density, so $G_{\nu}$ is essentially one, i.e., up to a constant.  The hypothesis of an existing continuous viscosity limit
\begin{equation}
g=\lim_{\nu\downarrow 0}g_{\nu}
\end{equation}
for $g_\nu\in S$ implies that the viscosity limit
\begin{equation}
\lim_{\nu\downarrow 0} g_{\nu}\ast G_{\nu}~\mbox{exists.}
\end{equation}
Indeed for $g_{\nu}\in S$ for $z_i=\frac{y_i}{\sqrt{\nu}},~1\leq i\leq D$ we have 
\begin{equation}\label{limnu}
\begin{array}{ll}
g_{\nu}\ast G_{\nu}=\int_0^t\int_{{\mathbb R}^D}g_{\nu}(s,x-y)\frac{1}{\sqrt{4\pi \nu s}^D}\exp\left(-\frac{|y|^2}{4\nu s} \right)dy ds=\\
\\
\int_0^t\int_{{\mathbb R}^D}g_{\nu}(s,x-\sqrt{\nu}z)\frac{1}{\sqrt{4\pi  s}^D}\exp\left(-\frac{|z|^2}{4 s} \right)dz ds\rightarrow \int_0^tg(s,x)ds \mbox{ as $\nu\downarrow 0$.}
\end{array}
\end{equation}
Similar considerations can be applied for multivariate spatial derivatives of order $\gamma$ if $D^{\gamma}_xg_{\nu}$ is at least continuous and $D^{\gamma}_xg=\lim_{\nu\downarrow 0}D^{\gamma}_xg_{\nu}$ exists.
In (\ref{Navlerayscheme}) for $k=1$ we have $v^{\nu,-,k-1}_i=v^{\nu,-,0}_i=v^{f,-}_i,~1\leq i\leq D$. Hence for $k=1$ the data $v^{f,-}_i$ in the convolutions $v^{f,-}_i\ast G_{\nu}$ and $v^{f,-}_i\ast_{sp} G_{\nu}$ are independent of $\nu$ and there are natural viscosity limits of these terms. For $k=2$ convolutions  of the form $D^{\gamma}_x\left( v^{f,-}_i\ast_{sp} G_{\nu}\right) $ appear in the recursive scheme (especially with $0\leq |\gamma|\leq 1$) and as part of the function $v^{\nu,-,1}_i$. We observe ( according to the previous remarks) that these functions $D^{\gamma}_x\left( v^{f,-}_i\ast_{sp} G_{\nu}\right)=\left( D^{\gamma}_x v^{f,-}_i\right) \ast_{sp} G_{\nu},~0\leq |\gamma|\leq 1$ have a $\nu$-independent upper bound. 

If we apply spatial derivatives to the scheme in (\ref{Navlerayscheme}), then it is natural to consider convolutions with first order spatial derivatives of the Gaussian $G_{\nu,i}$. We prove the existence of $\nu$-independent upper bounds $h_2,h_k,~k\geq 3$, where for some $T>0$
\begin{equation}
\sup_{t\in [0,T]}{\big |}\delta v^{\mbox{init},\nu,-,2}_i(t,.){\big |}_{H^3\cap C^3}\leq h_2,
\end{equation}
and where for $k\geq 3$
\begin{equation}
\sup_{t\in [0,T]}{\big |}(\delta v^{\nu,-,k}_i(t,.){\big |}_{H^3\cap C^3}\leq h_k.
\end{equation}
Here we shall use spatial Lipschitz continuity of the convoluted data $$\delta v^{\mbox{init},\nu,-,2}_i,\delta v^{\nu,-,k}_i,~1\leq i\leq  D, \mbox{for $k\geq 3$}.$$ Local time $\nu$-independent contraction in spatial $|.|_{H^3\cap C^3}$-norms of the higher order increments $\delta v^{\nu,-,k}_i,~1\leq i\leq  D$ for $k\geq 3$ leads to a regular limit for the functional series in (\ref{v-rep}) on the time interval $(0,T]$ for some time horizon $T>0$. Furthermore we shall observe that strong polynomial decay of order $2(D+1)$ is preserved for the increments, i.e., for $k\geq 3$ there is a $C>0$ independent of $\nu$ such that for all $|x|\geq 1$ and $0\leq |\gamma|\leq 3$
\begin{equation}
{\big |}D^{\gamma}_x\delta v^{\nu,-,k}_i(t,x){\big |}\leq \frac{C}{1+|x|^{2(D+1)}},
\end{equation}
and such that for all $|x|\geq 1$ and $0\leq |\beta|\leq 2$
\begin{equation}
{\big |}D^{\beta}_x\delta v^{\nu,-,k}_i(t,x){\big |}\leq \frac{C}{1+|x|^{2(D+1)}}.
\end{equation}
We can than use compactness arguments or $\nu$-independent local time contraction in order to verify the regularity of the limit for the functional series in (\ref{v-rep}).

In order to obtain $\nu$-independent estimates for convolutions of type $D^{\beta}_xg_{\nu}\ast G_{\nu,i}$ for $g_{\nu}\in S$ we may use  Lipschitz continuity of the convoluted data. For the first increment member $\delta v^{\mbox{init},\nu,-,2}_i$ in (\ref{v-rep}) this means
that for multiindices $0\leq |\beta|\leq 2$ and $\forall y,y'\in {\mathbb R}^D$~we have
\begin{equation}\label{c1}
~{\big |}D^{\beta}_x\delta v^{\mbox{init},\nu,-,2}_i(t,y)-D^{\beta}_x\delta v^{\mbox{init},\nu,-,2}_i(t,y'){\big |}\leq c{\big |}y-y'{\big |}
\end{equation}
for some finite constant $c>0$. 
 Lipschitz continuity turns out to hold  for the higher order increments $D^{\beta}_x\delta v^{,\nu,-,k}_i,~k\geq 3$ for $0\leq |\beta|\leq 2$ as well. 
We can have $c$ independent of $\nu$. Here we observe that the convolution with the Gaussian and the first order derivative of the Gaussian degenerates outside a ball $B_{\nu^{0,5-\epsilon}}(x)$ of radius $\nu^{0,5-\epsilon}$ for $0<\epsilon <0.5$ as $\nu$ becomes small.  
Indeed this holds for all spatial derivatives of the Gaussian as we have
\begin{equation}
\begin{array}{ll}
{\Big |}D^{\gamma}_xG_{\nu}(t,x;s,y){\Big |}{\Bigg |}_{|x-y|\geq \nu^{0.5-\epsilon}}\\
\\
\leq {\Big |}\frac{H_{\gamma,\nu (t-s)}(x-y)}{\sqrt{4\pi \nu(t-s)}^3} \exp\left(-\frac{|x-y|^2}{4\nu(t-s)} \right){\Big |}{\Bigg |}_{|x-y|\geq \nu^{0.5-\epsilon}}\downarrow 0~\mbox{as}~\nu\downarrow 0,
\end{array}
\end{equation}
where $H_{\gamma,\nu (t-s)}(x-y)$ is a Hermite-type polynomial of order $|\gamma|$ parameterized by $1/\nu (t-s)$.
If the viscosity limit $\lim_{\nu\downarrow 0}g_{\nu}\ast D^{\gamma}_xG_{\nu}(t,x)$ exists (is finite for some $x$) , then the 'mass' is concentrated in a ball of radius $B_{\nu^{0.5-\epsilon}}(x)$, i.e., assuming that the limits exist we have
\begin{equation}
\begin{array}{ll}
\lim_{\nu\downarrow 0}g_{\nu}\ast D^{\gamma}_{x}G_{\nu}(t,x)=\lim_{\nu\downarrow 0}\int_0^t\int_{{\mathbb R}^3}g_{\nu}(s,y)D^{\gamma}_{x}G_{\nu}(t,x;s,y)dyds\\
\\
=\lim_{\nu\downarrow 0}\int_0^t\int_{B_{\nu^{0.5-\epsilon}}(x)}g_{\nu}(s,y)D^{\gamma}_{x}G_{\nu}(t,x;s,y)dyds.
\end{array}
\end{equation}
Locally, we use simple Gaussian estimates of the form
\begin{equation}\label{gaussest}
\begin{array}{ll}
|D^{\gamma}G_{\nu}(t,x;s,y)|\\
\\
\leq \frac{{\big |}H_{\gamma,\nu (t-s)}(x-y){\big |}}{\sqrt{4\pi}^3}\frac{1}{\nu^{\delta} (t-s)^{\delta}}\left( |x-y
|^2\right)^{\delta -\frac{3}{2}} \left( \frac{|x-y
|^2}{\nu (t-s)}\right)^{\frac{3}{2}-\delta} \exp\left(-\frac{|x-y|^2}{4\nu(t-s)} \right)\\
\\
\leq \frac{C_{|\gamma|}}{\nu^{\delta}t^{\delta}|x-y|^{3+|\gamma|-2\delta}}
\end{array}
\end{equation}
for $\delta\in (0,1)$. We remark that for $0\leq |\gamma|\leq 1$ the constant $C_{|\gamma|}$ in (\ref{gaussest}) can be chosen to be   
\begin{equation}\label{c00}
C_0=\sup_{z>0, \delta \in (0,1)}z^{\frac{3}{2}-\delta}\exp\left(-\frac{z^2}{4}\right) 
\end{equation}
and 
\begin{equation}\label{c11}
C_1=\sup_{z>0, \delta \in (0,1)}z^{\frac{5}{2}-\delta}\exp\left(-\frac{z^2}{4}\right). 
\end{equation}
Note that $C_0,C_1$ in (\ref{c00}),~(\ref{c11}) are independent of $\nu$. 
For first order spatial derivatives and regular data $h$ with finite upper bound $C_h$ we observe that for any $\epsilon >0$
\begin{equation}\label{congkoma}
\begin{array}{ll}
{\Big |}\int_{y\in {\mathbb R}^3\setminus B^3_{\nu^{0.5-\epsilon}}(x)}h(y)\frac{(x-y)_j}{4\nu (t-s)\sqrt{4\pi \nu(t-s)}^3} \exp\left(-\frac{|x-y|^2}{4\nu(t-s)} \right)dy_1dy_2dy_3{\Big |}\\
\\
\leq {\Big |}\frac{cC_h}{\sqrt{4\pi \nu(t-s)}} \exp\left(-\frac{1}{4\nu^{2\epsilon}(t-s)} \right){\Big |}\downarrow 0 \mbox{ as }\nu\downarrow 0,
\end{array}
\end{equation} 
where we use Lipschitz continuity of the data, and where $c>0$ is a positive constant. The observation in (\ref{congkoma}) tells us that any mass of a convolution with $G_{\nu,j}$ evaluated at $(t,x)$ is concentrated in a ball $B_{\sqrt{\nu}}(x)$ of radius $\sqrt{\nu}$ around $x$ as $\nu$ becomes small.
In order to obtain $\nu$-independent estimates we may use Lipschitz continuity for terms 
\begin{equation}
 h\ast G_{\nu,k} \mbox{ (which are equivalent for $\nu>0$)}
\end{equation}
with data $h=D^{\beta}_xg_{\nu}$ for $g_{\nu}\in S$ and $0\leq |\beta| \leq 2$.
Observe the symmetry
\begin{equation}\label{symm}
G_{\nu,i}(t,x;s,y)=-G_{\nu,i}(t,x^{i,-};s,y^{i,-}),
\end{equation}
where $x=(x_1,\cdots ,x_n)$, $x^{i,-}=\left(x^{i,-}_1,\cdots,x^{i,-}_n\right)$ and $x^{i,-}_j=x_j-2\delta_{ij}x_j$ with the Kronecker $\delta_{ij}$ and where the symbol $y^{i,.}$ is defined analogously. Then Lipschitz continuity (or H\"{o}lder continuity with exponent $\delta_1=1$) of the shifted convoluted function $h_x(.)=h(x-.)$ and transformation $z=x-y$ (with $r=\sqrt{z_1^2+z_2^2+z_3^2}$ for brevity) leads to
\begin{equation}\label{visclim}
\begin{array}{ll}
{\Big |} \int_{B_{\nu^{0.5}}(x)}h(y)\frac{2(x_i-y_i)}{4\nu s\sqrt{4\pi \nu s}^D}\exp\left(-\frac{|x-y|^2}{4\nu s}\right) dy{\Big |}\\
 \\
\leq {\Big |}\int_{B_{\nu^{0.5}}, y_i\geq 0 }{\Big |}h_x(z)-h_x(z^{i,-}){\Big |}\frac{2z_j}{4\nu s\sqrt{\pi \nu s}^D}\exp\left(-\frac{|z|^2}{4\nu s}\right) dz{\Big |}\\
\\
\leq \int_{B_{\nu^{0.5}}, y_i\geq 0 }4|r|^{\delta_1}\frac{C_{1}}{\nu^{\delta}t^{\delta}|r|^{4-2\delta}}r^2dr\\
\\
4\frac{\tilde{C}}{\nu^{\delta}t^{\delta}}|r|^{\delta_1+2\delta-1}{\Big |}^{r=\nu^{0.5}}_0=4\frac{\tilde{C}}{t^{\delta}}|\nu|^{(0.5)(\delta_1+2\delta-1)-\delta}= 4\frac{\tilde{C}}{t^{\delta}}~\mbox{ $\forall~\nu \geq 0$ and $\delta_1=1$.}
\end{array}
 \end{equation}
Here, $\delta\in (0.5,1)$ and the latter equation holds even for $\nu=0$ if the usual definition $0^0=1$ is used.
This estimate can be used if the Leray data function, i.e., the Leray projection term applied to the data, is Lipschitz. Lipschitz Leray data function may correspond to a weak singularity of the vorticity, i.e. an oscillatory bounded vorticity which is not continuous. However, we shall observe that even H\"{o}lder continuous data for a velocity component can be chosen and there is still a global solution branch of the reversed Euler equation with such data. This global solution branch is obtained by a viscosity limit of a scheme derived from a Navier Stokes type equation.

%

For well chosen initial data $v^{f,-}_i,~1\leq i\leq 3$ (for the viscosity extension of the time reversed problem) we observe that $\delta v^{\mbox{init},\nu,-,2}_i$ and $\delta v^{\nu,-,2}_i$ gain regularity independently of $\nu$ as $k$ increases from $1$ to $2$ and observe that for all $t\geq 0$ the function
\begin{equation}\label{vseries}
v^{\nu,-}_i(t,.)-v^{f,-}_i\ast_{sp}G_{\nu}=\delta v^{\mbox{init},\nu,-,2}_i(t,.)+\sum_{l=3}^{\infty}\delta v^{\nu,-,l}_i(t,.)
\end{equation}
have $\nu$-independent spatial regular upper bounds with respect to the $H^3\cap C^3$-norm.
Moreover we shall observe that the functional series in (\ref{vseries}) can be differentiated twice with respect to the spatial variable member by member. Here, we may use strong polynomial decay at spatial infinity, which leads to stronger compactness arguments. The property of strong polynomial decay at spatial infinity is considered below (indeed, we shall observe that a spatial polynomial decay at spatial infinity of order $2(D+1)$ is preserved by the increments in (\ref{vseries}) and for spatial derivatives of these increments up to order $3$).  Moreover we shall observe that
\begin{equation}
\lim_{\nu\downarrow 0}\nu \Delta v^{\nu,-}_i\equiv 0,
\end{equation}
and conclude that $v^{\nu,-}_i,~1\leq i\leq D$ 
is a regular classical solution of the viscosity extension of the time-reversed incompressible Euler equation which has a regular viscosity limit which solves the time reversed incompressible Euler equation.

\begin{rem}
Note that the viscosity limit $\nu\downarrow 0$ for $D^{\beta}_xv^{\nu,-}_i,~1\leq i\leq D,~0\leq |\beta|\leq 2$ can be obtained by limits
\begin{equation}\label{vseries2}
D^{\beta}_xv^{\nu,-}_i-v^{f,-}_i\ast_{sp}D^{\beta}_xG_{\nu}=D^{\beta}_x\delta v^{\mbox{init},\nu,-,2}_i+\sum_{l=3}^{\infty}D^{\beta}_x\delta v^{\nu,-,l}_i,
\end{equation}
(evaluated at $t>0$), and where for $0\leq |\beta|\leq 2$ we can use for $t>0$
\begin{equation}\label{Navlerayscheme2v}
\begin{array}{ll}
D^{\beta}_x\delta v^{\mbox{init},\nu,-,2}_i= D^{\beta}_xv^{\nu,-,2}_i-v^{f,-}_i\ast_{sp}D^{\beta}_xG_{\nu}\\
\\
=\sum_{j=1}^D \left( D^{\beta}_{x}\left( v^{\nu,-,1}_j\frac{\partial v^{\nu,-,1}_i}{\partial x_j}\right) \right) \ast G_{\nu}\\
\\-\left( \sum_{j,m=1}^D\int_{{\mathbb R}^D}\left( D^{\beta}_x\left( \frac{\partial}{\partial x_i}K_D(.-y)\right) \left( \frac{\partial v^{\nu,-,1}_m}{\partial x_j}\frac{\partial v^{\nu,-,1}_j}{\partial x_m}\right)\right)  (.,y)dy\right) \ast G_{\nu},
\end{array}
\end{equation}
and variations of this formula obtained by the convolution rule. 
Furthermore, for $k\geq 3$ and $t>0$ we may use
\begin{equation}\label{Navlerayscheme2v3}
\begin{array}{ll}
D^{\beta}\delta v^{\nu,-,k}_i=
+\sum_{j=1}^D \left( D^{\beta}_x\left( v^{\nu,-,k-1}_j\frac{\partial v^{\nu,-,k-1}_i}{\partial x_j}\right)\right)  \ast G_{\nu}\\
\\-\left( \sum_{j,m=1}^D\int_{{\mathbb R}^D}\left( \frac{\partial}{\partial x_i}K_D(.-y)\right) \left( D^{\beta}_x\left( \frac{\partial v^{\nu,-,k-1}_m}{\partial x_j}\frac{\partial v^{\nu,-,k-1}_j}{\partial x_m}\right) (.,y)dy\right)\right)  \ast G_{\nu}\\
\\
-\sum_{j=1}^D \left( D^{\beta}_x\left( v^{\nu,-,k-2}_j\frac{\partial v^{\nu,-,k-2}_i}{\partial x_j}\right)\right)  \ast G_{\nu}\\
\\+\left( \sum_{j,m=1}^D\int_{{\mathbb R}^D}\left( D^{\beta}_x\left( \frac{\partial}{\partial x_i}K_D(.-y)\right) \left( \frac{\partial v^{\nu,-,k-2}_m}{\partial x_j}\frac{\partial v^{\nu,-,k-2}_j}{\partial x_m}\right) (.,y)dy\right)\right)  \ast G_{\nu},
\end{array}
\end{equation}
or variations by shifts of one spatial derivative to the Gaussian by the convolution rule.
For $0\leq |\beta|\leq 2$ we can consider the viscosity limit of a functional series of convolutions with the probability density $G_{\nu}$. There is no degeneracy issue for this classical limit.
\end{rem}

We summarize the six steps of the proof.
\begin{itemize}
\item[i)] As usual let $H^{2,1}$ be the Sobolev space where weak derivatives up to second order are in $L^1$. There are stronger singularities in this space in the space $H^2$. We provide examples of singular vorticity data in both spaces. We choose data $v^{f,-}_i,~1\leq i\leq D$, which determine the vorticity data $\omega^{f,-}_i=\omega^{\nu,-}_i(0,.),~1\leq i\leq D$ by the curl operation. Recall that the velocity data can be recovered from the vorticity data by the Biot-Savart law. First we cchoose data for one specific velocity component. In a second step we determine the data for the other velocity components according to different variations of argument. 

For a positive number $\beta_0\in (2,2+\alpha_0)$ with $\beta_0$ close to $2+\alpha_0$, and $\alpha_0\in \left(0,\frac{1}{2}\right)$ we consider  velocity component data $v^{f,-}_i\in H^{2,1},~1\leq i\leq 3$, where for one index $i_0\in \left\lbrace 1,2,3\right\rbrace$  we have $v^{f,-}_{i_0}(x)=g_{(0)}(r)$ for some univariate function $g_{(0)}$ which we are going to define next. Here the function $g_{(0)}$ is dependent on the radial component of spherical polar coordinates $r=\sqrt{x_1^2+x_2^2+x_3^2}\in {\mathbb R}_+$, where ${\mathbb R}_+$ denotes the set of nonnegative real numbers. This function  $g_{(0)}:{\mathbb R}_+\rightarrow {\mathbb R}$ is determined by
\begin{equation}\label{g0}
g_{(0)}(r)= \phi_{1}(r)r^{\beta_0}\sin\left(\frac{1}{r^{1+\alpha_0}}\right),~g_{(0)}(0)=0,
\end{equation}
where $\phi_1 \in C^{\infty}_0$ (the space of smooth functions wanshing at infinity) may be defined by
\begin{equation}
\phi_{1}(r)=\left\lbrace \begin{array}{ll}
1~\hspace{1.75cm}\mbox{if }~r\leq 1,\\
\phi_1(r)=\alpha_*(r) \mbox{ if }~1\leq r\leq 2,
\\
0~\hspace{1.65cm}\mbox{ if }~r\geq 2.
\end{array}\right.
\end{equation}
Here, $\alpha_*$ is a smooth function with bounded derivatives for $1\leq r\leq 2$. This choice of $\phi_1$ is in the scape $C^{\infty}_c$ of smooth functions with compact support.

\begin{rem}
Instead of $\phi_1$ the function 
\begin{equation}\label{phi}
\phi(r)=\frac{1}{(1+r^2)^2},~ r=\sqrt{x_1^2+x_2^2+x_3^3}
\end{equation}
is an alternative choice which works for all variations of arguments considered here.
\end{rem}
\begin{rem}
An alternative function  $g^*_{(0)}:{\mathbb R}_+\rightarrow {\mathbb R}$ determined by
\begin{equation}
g^*_{(0)}(r)= \phi_{1}(r)r^{\beta^*_0}\sin\left(\frac{1}{r^{\alpha_0}}\right),~g_{(0)}(0)=0,
\end{equation}
where $\beta^*_0\in (1,1.5)$ and $\phi_1 \in C^{\infty}_0$ is as before or, alternatively replaced by (\ref{phi}), defines a function $g^*_{0}\in H^2$ with a singular vorticity at the origin. This singularity is slightly weaker than the singularity of the function $g_{(0)}\in H^{2,1}$, but it suffices for the purposes of this article as well.  
\end{rem}
Functions as $\phi_1$ are well-known in the context of partitions of unity. As usual, $C^{\infty}_c$ denotes the function space of smooth functions with compact support. In general $\beta_0-1$ will be chosen close to $1+\alpha_0$. Note that the data function $v^{f,-}_{i_0}$ is only locally Lipschitz. The regularity gain of the local solution $v^{\nu,-}_i,~ 1\leq i\leq D$ is obtained for the solution increment $\delta v^{\nu,-}_i=v^{\nu,-}_i-v^{f,-}_i$.

\begin{rem}
In case $\beta_0=2+\alpha_0$ we have Lipschitz continuous data. The singularities of the vorticities are then bounded and of oscillatory type. The following argument holds also in this case, but the result is weaker. 
\end{rem}

\begin{rem}
In order to prove the existence of kinks or weak singularities we may consider for $k\geq 2$ and $\beta^k_0\in (k+1,k+1+\alpha_0)$, $\alpha_0\in \left(0,\frac{1}{2}\right)$ velocity component data $v^{f,-}_i\in H^{2+k},~1\leq i\leq 3$ where for one index $i_0\in \left\lbrace 1,2,3\right\rbrace$  we have $v^{f,-}_{i_0}(x)=g_{(k)}(r)$ and where $g_{(k)}:{\mathbb R}_+\rightarrow {\mathbb R}$ is defined by
\begin{equation}
g_{(k)}(r)= \phi_{1}(r)r^{\beta^k_0}\sin\left(\frac{1}{r^{1+\alpha_0}}\right).
\end{equation}
\end{rem}

Now we have constructed data for one velocity component $v^{f,-}_{i_0}$ with a singularity which becomes a final data component for the inverse Euler equation. Recall that the data term in the local solution representation of the related Navier Stokes type extension of the time-revesed Eulerequation is of the form $v^{f,-}_{i_0}\ast_{sp}G_{\nu}$. A viscosity limit $\nu_k\downarrow 0$ yields that data $v^{f,-}_{i_0}$ for this data term part of a local solution for the time reversed Euler equation. However, we have observed that the corresponding term first order spatial derivatives data term $v^{f,-}_{i_0}\ast_{sp}G_{\nu,i_0}$ degenerates (as we have a spatial convolution). This leads to several possibilities concerning the choice of the data of the other velocity components $v^{f,-}_j,~j\neq i_0$.
\begin{rem}\label{datarem}
We can use data of similar regularity, i.e., the construction of singularities can also be done with the data 
\begin{equation}
v^{f,-}_1(x) = x_2x_3
\frac{r^{\beta^*_0}\sin\left(r^{-\alpha_0}\right)}{(1+r^2)^2},~
v^{f,-}_2(x) =-\frac{1}{2} x_1x_3
\frac{r^{\beta^*_0}\sin\left(r^{-\alpha_0}\right)}{(1+r^2)^2}
\end{equation}
and
\begin{equation}
v^{f,-}_3(x)=-\frac{1}{2}x_1x_2
\frac{r^{\beta^*_0}\sin\left(r^{-\alpha_0}\right)}{(1+r^2)^2},~ r=\sqrt{x_1^2+x_2^2+x_3^3}.
\end{equation}
Here $\beta^*_0\in (1,1.5)$ and $\alpha_0\in (0,0.5)$.
\end{rem}
Next we consider choices of data related to the choice (\ref{g0}) above.
For $j\in \left\lbrace 1,2,3\right\rbrace\setminus \left\lbrace i_0 \right\rbrace$ we may choose velocity component data such that
\begin{equation}\label{incompresscond}
\sum_{j=1}^Dv^{f,-}_{j,j}=0,~~\sum_{j=1}^D\omega^{f,-}_{j,j}=0
\end{equation}
is satisfied in a distributive sense (and in a classical sense outside the origin). We observed above that for a local solution increment $\delta v^{\nu,-}_i,~ 1\leq i\leq D$ which is regular (spatially in $H^2\cap C^2$ and independently of $\nu$),  the incompressibility condition holds in the viscosity limit $\nu_k\downarrow 0$ over a time intevral if we choose functions $v^{f,-}_j,~ 1\leq j\leq D,~ j\neq i_0$ such that for all  $t>0$
\begin{equation}
\sum_{i=1}^Dv^{f,-}_k\ast_{sp}G_{\nu_k,i}(t,.)\downarrow 0,~\mbox{as}~ \nu_k\downarrow 0.
\end{equation}
This observation implies that we have some freedom in the choice of data even if $v^{f,-}_{i_0}$ is chosen. We consider a small list
\begin{itemize}
\item[ia)] we can choose data as in (\ref{datarem}). These are $H^2$ data where all velocity components are of similar regularity. The singularity of the vorticity is at one point in space time. However all vorticity components have a singularity. In this case the argument below  does not depend on the degeneracy of the first order spatial derivatives of the convoluted initial data.

\item[ib)] Let $i_0=1$and let $v^{f,-}_{i_0}=v^{f,-}_1$ be chosen as in (\ref{g0}). Siince $v^{f,-}_{i_0}$ is locally regular around any argument $x\neq 0$, we have the degeneracy  $v^{f,-}_{i_0}\ast_{sp}G_{\nu_k,i}(t,x)\downarrow 0,~\mbox{as}~ \nu_k\downarrow 0$ for $x\neq 0$ and also a degeneracy at $x=0$ by the symmetry of the data. Hence, if we have regularity of the local solution increment $\delta v^{\nu_k,-}_i,~ 1\leq i\leq D$ (a regularity which is spatiall at least $H^2\cap C^2$), which holds also in the viscosity limit $\nu_k\downarrow 0$, then  the condition in (\ref{incompresscond}) reduces to
$\sum_{j=1,j\neq 1}^Dv^{f,-}_{j,j}=0$. Then we may choose regular data $v^{f,-}_j,~ 2\leq j\leq 3$ such that ib$\alpha$) the latter reduced incompressibility condition is satisfied, and ib$\beta$) the local time vorticity solution increment does not degenerate in the viscosity limit. 
For example we may choose $v^{f,-}_2(x):=x_3r^{\beta_0}\phi_1$ and $v^{f,-}_3(x)=-x_2r^{\beta_0}\phi_1$ where $\phi_1$ may be chosen as above and such that it depends only on the radial variable $r=\sqrt{x_1^2+x_2^2+x_3^2}$. In this case we have a stronly regularnontrivial local solution increment for the velocity and for the vorticity, where the incompressibility condition and the dimension $D=3$ is explicityly used. The regularity argument for the solution increment can then be based on an alternative argument. Note that for the choice in item ib) the incompressibility condition is only satosfied in the viscosity limit $\nu_k\downarrow 0$, and approximately satisfied for small $\nu_k>0$. However, this is sufficient for our purposes as we are interested in the viscosity limit. Furthermore note that the specific choice in this item takes advantage of the fact that the dimension is $D=3$ (or $D\geq 3$, where similar constructions are possible).

\item[ic)]Note that with our choice in this item we  have for any $\nu >0$ and any $t>0$
\begin{equation}
v^{f,-}_{i_0}\in H^1,~\mbox{therefore}~v^{f,-}_{i_0}\ast_{sp}G_{\nu,i_0}(t,.)=v^{f,-}_{i_0,i_0}\ast_{sp}G_{\nu}(t,.)
\end{equation}
Therefore, we may choose any $v^{f,-}_j,~j\neq i_0$ such that for any $\nu >0$ and any $t>0$
\begin{equation}
v^{f,-}_{i_0,i_0}\ast_{sp}G_{\nu}(t,.)=-\sum_{j\neq i_0}v^{f,-}_{j,j}(t,.)\ast_{sp}G_{\nu}(t,.)
\end{equation}
An example of such an approriate choice is to choose $v^{f,-}_{1}=v^{f,-}_{i_0}$ as above, and
\begin{equation}
v^{f,-}_2\ast_{sp}G_{\nu}(t,.)=-x_2v^{f,-}_{1,1}\ast_{sp}G_{\nu}
\end{equation}
for the second velocity component. Then we have
\begin{equation}
v^{f,-}_{1,1}\ast G_{\nu}+v^{f,-}_{2,2}\ast G_{\nu}=-x_2v^{f,-}_{1,1}\ast_{sp}G_{\nu,2}\stackrel{!}{=}v^{f,-}_{3,3}\ast_{sp}G_{\nu}
\end{equation}
The latter equation is an integral equation for $v^f_3$ which can be solved (maybe even explicitly for appropriate choices of $\phi$).
\end{itemize}

We observe that for our choices we have an upper bound in in $H^3\cap C^3$ which is independent of $\nu$.
 
We have already observed that the local solution of a Navier Stokes type equation (related to the time reversed Euler equation) satisfies incompressibility for $t>0$ (as does the viscosity limit). 
Note that we cannot impose a incompressible condition at $t=0$ in a classical sense, as we have a singular point at time $t=0$.
    
Note that for $k=0$ and $\beta_0=\beta^0_0$ close to $1+\alpha_0$ we have $v^{f}_i\in H^{2,1}$ $\omega^f_i\in H^{1,1}$ for all $1\leq i\leq D=3$ as the singularity of $v^{f}_{i_0}=v^{f}_{i_0}$ (resp. $\omega^f_{i_0}=\omega^f_1$) at $r=0$ (polar coordinates) is of order $|\beta^0_0-4-2\alpha_0|<2.5$ for $\beta_0$ close to $2+\alpha_0$ as $\alpha_0\in (0,0.5)$. Similarly for any $l\geq 0$ we have choices $\beta^k_0\in (k+1,k+1+\alpha]_0$, $\alpha_0\in \left(0,\frac{1}{2}\right)$ such that 
\begin{equation}
g_{(k)}(.)\in C^{l-1}\setminus C^{l}.
\end{equation}
We just have to choose $l=\min \left\lbrace m\geq 0| \beta^k_0-m(2+\alpha_0)<0\right\rbrace $. For this choice of $l$ we have $g_{(k)}(.)\in H^{l,1}$.

\item[ii)] 
As pointed out in the previous item, the choice of the data depends on the regularity of a nontrivial local solution increment $\delta v^{\nu_k,-}_i,~1\leq i\leq D$ where regularity and nontriviality hold in the viscosity limit $\nu_k\downarrow 0$. The argument for the regularity of this increment simplifies or can be stengthened if we choose the data as in ib$\beta$) in item i) above. However in any case, i.e., for any choice of data as listed in item i) we can observe the following.  
For $D\geq 3$ and for given $\nu >0$, $1\leq i\leq 3$ and for all , $0\leq |\gamma|\leq 1$, $x\in {\mathbb R}^3\setminus \left\lbrace 0\right\rbrace$,~$r=\sqrt{x_1^2+x_2^2+x_3^2}\in {\mathbb R}\setminus \left\lbrace 0\right\rbrace$  we have for the first iteration ($k=1$) 
\begin{equation}\label{vinit0}
~{\big |}D^{\gamma}_x\delta v^{\mbox{init},\nu,-,1}_i(t,x){\big |}={\big |}D^{\gamma}_x\delta v^{\nu,-,1}_i(t,x){\big |}\leq Cr^{2\beta_0-1-|\gamma|},
\end{equation}
and  for the second iteration ($k=2$) and  all $\delta \in (0.5,1)$
\begin{equation}\label{vinit1}
~{\big |}D^{\gamma}_x\delta v^{\mbox{init},\nu,-,2}_i(t,x){\big |}\leq Cr^{2\beta_0-1-|\gamma|}.
\end{equation}
For the higher order increments and for $k\geq 3$ we also have for  all $\delta \in (0.5,1)$
\begin{equation}
~{\big |}D^{\gamma}_x\delta v^{\nu,-,k}_i(t,x){\big |}\leq Cr^{2(\beta_0-1)+1-|\gamma|}.
\end{equation}
Indeed these increments gain regularity at each iteration step, but we shall not consider this here.

Inheritance of polynomial decay (cf. item iii)) then implies for $1\leq i\leq 3$ and for all $k\geq 3$ 
\begin{equation}
\delta v^{\nu,-,k}_i(t,.)\in  H^3\cap C^3,
\end{equation}
and
\begin{equation}
\delta v^{\mbox{init},\nu,-,2}_i(t,.)\in  H^3\cap C^3.
\end{equation}
Note that at this stage we do not consider dependence on $\nu$. The independence of $\nu$ requires an additional argument and is considered in item v).
\item[iii)] For a short time horizon $T>0$ for all $\nu >0$ polynomial decay of order  $2( D+1)$ at spatial infinity is inherited by the  increments $$D^{\gamma}_x\delta v^{\mbox{init},\nu,-,2}_i,D^{\gamma}_x\delta v^{\nu,-,k}_i,~k\geq 3,~0\leq |\gamma|\leq 3.$$ More precisely for $1\leq i\leq 3$, and $k\geq 3$ there exists a finite constant $c>0$  such that for all $t\in [0,T]$ and $|x|\geq 1$  
\begin{equation}\label{polybound}
\forall~ k\geq 1~\forall 0\leq|\gamma|\leq 3~{\Big |}D^{\gamma}_x\delta v^{\mbox{init},\nu,-,2}_i(t,x){\Big |}\leq  \frac{c}{1+|x|^{2(D+1)}},
\end{equation} 
and
\begin{equation}\label{polybound2}
\forall~ k\geq 1~\forall 0\leq|\gamma|\leq 3~{\Big |}D^{\gamma}_x\delta v^{\nu,-,k}_i(t,x){\Big |}\leq  \frac{c}{1+|x|^{2(D+1)}}.
\end{equation} 
 For small $T>0$ it is observed that the constant $c>0$ can be chosen independently of the iteration index $k$. Furthermore, for $0\leq |\gamma|\leq 4$ 
\begin{equation}
{\Big |}v^{f,-}_i\ast_{sp} D^{\gamma}_xG_{\nu}(t,.){\Big |}\leq  \frac{c}{1+|x|^{2(D+1)}}~\mbox{$t>0$ and $c$ depends on $t>0$},
\end{equation}
such that a corresponding statement holds for the value functions $v^{\nu,-,k}_i,~1\leq i\leq 3$ and $k\geq 2$.

\item[iv)] For the parameter choices of $\beta_0,\alpha_0$ with $\beta^*_0=1.5-\epsilon$ for small $\epsilon >0$ of item i), for $1\leq i\leq 3$, $\nu >0$ and for $k\geq 3$ and some time $T>0$
\begin{equation}\label{contr}
{\big |}\delta v^{\nu,-,k+1}_i(t,.){\big |}_{H^{3}\cap C^3}\leq \frac{1}{2}{\big |}\delta v^{\nu,-,k}_i(t,.){\big |}_{H^{3}\cap C^3}.
\end{equation}
Here we use classical representations of solution increments in form of convolutions of the Burgeres term and the Lery projection term with first order spatial derivatives of the Gaussian. Here we use Lipschitz continuity of the Burgers term and the Leray projection term established by an analyis via local time iteration schemes. Similarly for spatial derivatives.  We underline and emphasize that there is {\it no degeneracy for the solution increment} $\delta v^{\nu,-}_i,~ 1\leq i\leq D$ {\it for any choice of data listed in item i)}, although there is degeneracy of the spatial first order derivatives of the convoluted data terms (as we have observed above).      
Hence, for all $t\in [0,T]$
\begin{equation}
v^{\nu,-}_i(t,.)-v^{f,-}_i\ast_{sp}G_{\nu}(t,.)=\delta v^{\mbox{init},\nu,-,2}_i+\sum_{l=3}^{\infty}\delta v^{\nu,-,l}_i\in H^3\cap C^3.
\end{equation}
Similarly for the choice $\beta_0$ instead of $\beta^*_0$ (cf. item i). Independence of $\nu$ is a further step, but Lipshitz continuity of (spatial derivatives of some orderof ) the Leray projection term and the Burgers term in local time representations of solution increments in terms of first order spatial derivatives is crucial, too. It is  considered in the next step.
\item[v)] There exists a $T>0$ independent of $\nu>0$ and a finite constant $C>0$ independent of $\nu >0$ such that 
\begin{equation}
\sup_{t\in (0,T]}{\big |}v^{\nu,-}_i(t,.)-v^{f,-}_i\ast_{sp}G_{\nu}(t,.){\big |}_{H^3\cap C^3}\leq C.
\end{equation}
Furthermore, for such $T>0$ independent of $\nu >0$ such that 
\begin{equation}\label{incompnote}
\mbox{ for all }~t\in (0,T]:~\sum_{i=1}^Dv^{\nu,-}_{i,i}(t,.)=0.
\end{equation}
We note that the relation in (\ref{incompnote}) hold for all $\nu >0$. The trick here is to use the smooth coefficients in the iteration scheme in order to estimate the functional increments (cf. below ad vi)).
\item[vi)] Choose a time horizon $T>0$ as in the previous step such that contraction holds for the higher order increments $\delta v^{\nu,-,k}_i$ with $k\geq 3$ as in (\ref{contr}). 
For  an upper bound $C>0$ independent of $\nu>0$ we have
\begin{equation}
\begin{array}{ll}
{\big |}\delta v^{\mbox{init},\nu,-,2}_i(t,.)+\sum_{l=3}^{\infty}\delta v^{\nu,-,l}_i(t,.){\big |}_{H^3\cap C^3}\\
\\
={\big |}v^{\nu,-}_i(t,.)-v^{f,-}_i\ast_{sp}G_{\nu}(t,.){\big |}_{H^3\cap C^3}\leq C.
\end{array}
\end{equation}
For all $\nu >0$ and for some $T>0$ (indpendent of $\nu$) the function
\begin{equation}
v^{\nu,-}_i(t,.)=v^{f,-}_i\ast_{sp}G_{\nu}(t,.)+\delta v^{\mbox{init},\nu,-,2}_i+\sum_{l=3}^{\infty}\delta v^{\nu,-,l}_i
\end{equation}
satisfies the Navier stokes equation on the time interval $[0,T]$ such that
\begin{equation}\label{Navlerayeq2**}
\begin{array}{ll}
\nu\Delta v^{\nu,-}_i= \frac{\partial v^{\nu,-}_i}{\partial t}
-\sum_{j=1}^D v^{\nu,-}_j\frac{\partial v^{\nu,-}_i}{\partial x_j}\\
\\+\sum_{j,m=1}^D\int_{{\mathbb R}^D}\left( \frac{\partial}{\partial x_i}K_D(.-y)\right) \sum_{j,m=1}^D\left( \frac{\partial v^{\nu,-}_m}{\partial x_j}\frac{\partial v^{\nu,-}_j}{\partial x_m}\right) (t,y)dy,
\end{array}
\end{equation}
where the right side of (\ref{Navlerayeq2**}) is the time-reversed Euler equation operator applied to $v^{\nu,-}_i$. We choose a sequence $(\nu_p)_{p\geq 1}$ converging to zero and use probability density representations (convolutions with $G_{\nu}$ itself) of increments $D^{\beta}_x\delta v^{\mbox{init},\nu_p,-,2}_i$ and $D^{\beta}_x\delta v^{\nu_l,-,l}_i,~l\geq 3$ for $0\leq |\beta|\leq 2$ in
\begin{equation}
v^{-,\nu_p}_i:=v^{f,-}_{i}\ast_{sp}G_{\nu}+\delta v^{-\nu_p}_i
\end{equation}
in order to show that the viscosity limit $v^{-}_i,~1\leq i\leq D$ is a classical solution of the time-reversed Euler equation.exists . Here
\begin{equation}
\delta v^{-}_i:=\lim_{\nu_p\downarrow 0 }\delta v^{\mbox{init},\nu_p,-,2}_i+\sum_{l=3}^{\infty}\delta v^{\nu_l,-,l}_i\in C^1\left(\left(0,T\right] , H^3\cap C^3\right). 
\end{equation} 
Furthermore the function $v^{-}_i=v^{f,-}_i+\delta v^{-}_i,~1\leq i\leq D$ solves the time-reversed incompressible Euler equation and satisfies the incompressibility condition in the time interval $(0,T]$, i.e.,
\begin{equation}
\mbox{ for all $t\in (0,T]$}~\sum_{i=1}^D v^{-}_{i,i}(t,.)=0.
\end{equation}

We conclude that the original Euler equation develops in opposite time direction from data at time $\tau=0$ (corresponding to time $s=T$ of the reversed-time Euler equation) a weak singularity at time $T>0$ (corresponding to data at time $t=0$ the time-reversed Euler equation). 
\end{itemize}

ad i) first we note that multiindices are denoted by $\alpha=(\alpha_1,\cdots,\alpha_D),\beta=(\beta_1,\cdots,\beta_D),\cdots$ in this paper, and $\alpha_0,\beta_0$ just refer to positive real numbers. For the derivative of the data $v^{f,-}_{i_0}=v^{f}_1$ and $k=0$ we compute for $r\neq 0$ and $r\leq 1$
\begin{equation}
\begin{array}{ll}
g'(r)=\frac{d}{dr} r^{\beta_0}\sin\left(\frac{1}{r^{1+\alpha_0}}\right)
=\beta_0 r^{\beta_0-1}\sin\left(\frac{1}{r^{1+\alpha_0}}\right)-(1+\alpha_0) r^{\beta_0-2-\alpha_0}\cos\left(\frac{1}{r^{1+\alpha_0}}\right).
\end{array}
\end{equation}
The derivative $g'$ of the function $g$ at $r=0$ is strongly singular for $\beta_0\in (2,2+\alpha_0)$ and $\alpha_0\in \left(0,\frac{1}{2}\right)$. Note that it is 'oscillatory' singular bounded for $\beta_0 =\alpha_0 \in \left(0,\frac{1}{2}\right)$. Note that for data $v^{f}_{i_0}(x_1,x_2,x_3)=g(r)$ we have (for $r\neq 0$)
\begin{equation}
v^{f,-}_{i_0,j}=g'(r)\frac{\partial r}{\partial x_j}=g'(r)\frac{x_j}{r}.
\end{equation}
In polar coordinates $(r,\theta,\phi)\in [0,\infty)\times [0,\pi]\times[0,2\pi]$ with 
\begin{equation}
x_1=r\sin(\theta)\cos(\phi),~x_2=r\sin(\theta)\sin(\phi),~x_3=r\cos(\theta),
\end{equation}
(where for $r\neq 0$ and $x_1\neq 0$ we have $r=\sqrt{x_1^2+x_2^2+x_3^2}, ~\theta=\arccos\left(\frac{x_3}{r}\right),~\phi=\arctan\left(\frac{x_2}{x_1} \right)$) we get
\begin{equation}
\begin{array}{ll}
v^{f,-}_{i_0,1}=g'(r)\frac{x_1}{r}=g'(r)\sin(\theta)\cos(\phi),\\
v^{f,-}_{i_0,2}=g'(r)\frac{x_2}{r}=g'(r)\sin(\theta)\sin(\phi),\\
v^{f,-}_{i_0,3}=g'(r)\frac{x_3}{r}=g'(r)\cos(\theta),
\end{array}
\end{equation}
such that we have
\begin{equation}
v^{f,-}_{i_0}\in H^{1,1} ~\mbox{obviously.}
\end{equation}
The second derivative of $g$ is
\begin{equation}
\begin{array}{ll}
g''(r)=\frac{d^2}{dr^2} r^{\beta_0}\sin\left(\frac{1}{r^{1+\alpha_0}}\right)\\
\\
=\frac{d}{dr}\left( \beta_0 r^{\beta_0-1}\sin\left(\frac{1}{r^{1+\alpha_0}}\right)-(1+\alpha_0) r^{\beta_0-2-\alpha_0}\cos\left(\frac{1}{r^{1+\alpha_0}}\right)\right)\\
\\
=\beta_0(\beta_0-1)r^{\beta_0-2}\sin\left(\frac{1}{r^{1+\alpha_0}}\right)
-(1+\alpha_0)\beta_0 r^{\beta_0-3-\alpha_0}\cos\left(\frac{1}{r^{1+\alpha_0}}\right)\\
\\
+(1+\alpha_0)(2+\alpha_0-\beta_0) r^{\beta_0-2-\alpha_0}\cos\left(\frac{1}{r^{1+\alpha_0}}\right)\\
\\
-(1+\alpha_0)^2 r^{\beta_0-4-2\alpha_0}\sin\left(\frac{1}{r^{1+\alpha_0}}\right).
\end{array}
\end{equation}
We have $v^{f,-}_{i_0}\in H^{2,1}$, since
\begin{equation}
\beta_0-4-2\alpha_0 >-3.
\end{equation}
Note that
\begin{equation}
v^{f,-}_{i_0}(0,.) \in C^{\delta}\left({\mathbb R}^3\right), 
\end{equation}
for H\"{o}lder constants of order $\delta \in \left(0,\frac{1}{2}\right)$.
Similarly we have
\begin{equation}
\omega^{f,-}_{i_0}=\omega^{f,-}_1=\omega_{1}(0,.)\in H^{1,1}.
\end{equation}
Note that in item i) and in subitem ia) of item i) we have listed some possible variations. The data $v^{f,-}_j$ for $j\neq i_0$ can be chosen chosen form the list of subitems ia), ib) or ic) in item i) above. Note that for the choice in subitem ib) the incompressibility condition is only satisfied in the viscosity limit (but this suffices for our purposes). 

ad ii) First we consider estimates which hold for all choices of data in items ia), ib), and ic) above. The choice of data in item ib) allows for a stronger variation of argument in the viscosity limit. For $k=1$ (\ref{Navlerayscheme}) becomes
\begin{equation}\label{Navlerayscheme1}
\begin{array}{ll}
 v^{\nu,-,1}_i=v^{f,-}_i(.)\ast_{sp}G_{\nu}
+\sum_{j=1}^D \left( \left( v^{f,-}_j\ast_{sp} G_{\nu}\right) \left(  v^{f,-}_i\ast_{sp}G_{\nu,j}\right) 
\right) \ast G_{\nu}\\
\\-\sum_{j,m=1}^D\int_{{\mathbb R}^D}\left( \frac{\partial}{\partial x_i}K_D(.-y)\right)\times\\
\\
\times \sum_{j,m=1}^D\left(\left(   v^{f,-}_m\ast_{sp}G_{\nu,j}\right)  \left( v^{f,-}_j\ast_{sp}G_{\nu,m}\right) \right) (t,y)dy\ast G_{\nu}\\
\\
=:v^{f,-}_i(.)\ast_{sp}G_{\nu}
+B_0\ast G_{\nu}-L_0\ast G_{\nu},
\end{array}
\end{equation}
where $B_0$ and $L_0$ denote abbreviations of the lowest order approximations of the Burgers term and the Leray projection term. Recall that at this stage we consider a given $\nu>0$, and do not consider dependence of $\nu$. Nevertheless the estimates are of a form such that an additional argument in item v) leads to $\nu$-independence of upper bounds. For $t\geq 0$ we have
\begin{equation}
\mbox{for all}~r=\sqrt{x_1^2+x_2 ^2+x_3^2}~\mbox{we have}~{\Big |} \left( v^{f,-}_j\ast_{sp} G_{\nu}\right) (t,r){\Big |}\leq cr^{\beta_0},
\end{equation}
and
\begin{equation}\label{vfi}
\mbox{for all}~r=\sqrt{x_1^2+x_2 ^2+x_3^2}~\mbox{we have}~{\Big |} \left(  v^{f,-}_i\ast_{sp}G_{\nu,j}\right)(t,.){\Big |}\leq cr^{\beta_0-1}
\end{equation}
for some finite constant $c>0$, which is independent of $\nu >0$. 
Consequently, we have for all $r=\sqrt{x_1^2+x_2 ^2+x_3^2}$
\begin{equation}
{\Big |} \left( v^{f,-}_j\ast_{sp} G_{\nu}\right) \left(  v^{f,-}_i\ast_{sp}G_{\nu,j}\right)(t,x){\Big |}\leq cr^{2\beta_0-1}
\end{equation}
and as $\frac{\partial}{\partial x_i}K_{D,i}(.-y)\sim \frac{1}{r^2}$ for $D=3$ we have for all $r=\sqrt{x_1^2+x_2 ^2+x_3^2}$
\begin{equation}
\begin{array}{ll}
{\Big |} \int_{{\mathbb R}^D}K_{D,i}(x-y) \sum_{j,m=1}^D\left( \left(  v^{f,-}_j\ast_{sp}G_{\nu,m}\right)\left(  v^{f,-}_m\ast_{sp}G_{\nu,j}\right)\right)  (t,y)dy{\Big |}\\
\\
\leq cr^{2\beta_0-2+1}=cr^{2\beta_0-1}.
\end{array}
\end{equation}

Using Gaussian estimates above for multiindices $\gamma$ with $0\leq |\gamma|\leq 1$
for the convoluted Burgers term $B^0$ we get 
\begin{equation}
{\Big |}B^0\ast D^{\gamma}_xG_{\nu}(t,.){\Big |}\leq cr^{2\beta_0-1-|\gamma|},
\end{equation}
and as $\frac{\partial}{\partial x_i}K_{D,i}(.-y)\sim \frac{1}{r^2}$ for $D=3$ for the convoluted Leray projection term $L^0$ we have
\begin{equation}
{\Big |}L^0\ast D^{\gamma}_xG_{\nu}(t,.){\Big |}\leq cr^{2\beta_0-1-|\gamma|}
\end{equation}
Hence the statement in (\ref{vinit0}) follows. Next we observe thatfor all $1\leq i\leq D$ and
\begin{equation}
\forall t>0~ ~\delta v^{\nu,-,1}_i(t,.)= v^{\nu,-,1}_i(t,.)-v^{f,-}_i\in C^2\cap H^2.
\end{equation}
First we consider the Leray projection term. For the second order spatial derivatives with indices $p,q$ we have 
\begin{equation}\label{leray}
\begin{array}{ll}
L^*_{i,p,q}:=L^0_{i,p}\ast G_{\nu,q}=-\sum_{j,m=1}^D\int_{{\mathbb R}^D}K_{D,i,p}(.-y)\times\\
\\
\times \sum_{j,m=1}^D\left(\left(   v^{f,-}_m\ast_{sp}G_{\nu,j}\right)  \left( v^{f,-}_j\ast_{sp}G_{\nu,m}\right) \right) (t,y)dy\ast G_{\nu,q}
\end{array}
\end{equation}
The singularity order of $K_{D,p,q}$ at $r=0$ is $r^{-n}$ such that for all $x$ in a finite ball around zero (recall $r=|x|=\sqrt{\sum_{j=1}^Dx_j^2}$) we have for all $t>0$ (fixed)
\begin{equation}\label{leraybb}
\begin{array}{ll}
{\Big |}\sum_{j,m=1}^D\int_{{\mathbb R}^D}K_{D,i,p}(x-y) \sum_{j,m=1}^D\left(\left(   v^{f,-}_m\ast_{sp}G_{\nu,j}\right)  \left( v^{f,-}_j\ast_{sp}G_{\nu,m}\right) \right) (t,y)dy{\Big |}\\
\\
\leq cr^{2\beta_0-2}
\end{array}
\end{equation}
where $c$ is independent of $\nu >0$ (and of $t$). 
Hence, for all $x$ in a finite ball around the origin we have 
\begin{equation}
{\Big |}L^0_{i,p,q}(t,x){\big |}\leq cr^{2\beta_0-3}.
\end{equation}
 The latter estimate holds a fortiori for lower order derivatives of $L_i$. It follows that  $L_i$ and spatial derivatives of $L_i$ up to second order are Lipschitz, especially continuous. Here, note that we easily conclude that for all $t>0$
\begin{equation}
L^*_{i,p,q}(t,.)\in L^2.
\end{equation}
This holds a fortiori for the lower order spatial derivatives such that we conclude that for all $t>0$
\begin{equation}
L^*_i(t,.)\in H^2\cap C^2.
\end{equation}
Next consider the Burgers term
\begin{equation}
B^*_i=\sum_{j=1}^D \left( \left( v^{f,-}_j\ast_{sp} G_{\nu}\right) \left(  v^{f,-}_i\ast_{sp}G_{\nu ,j}\right) 
\right) \ast G_{\nu}.
\end{equation}
We have
\begin{equation}\label{bb*}
\begin{array}{ll}
B^*_{i,p,q}=\sum_{j=1}^D \left( \left( v^{f,-}_j\ast_{sp} G_{\nu,p}\right) \left(  v^{f,-}_i\ast_{sp}G_{\nu,j}\right) 
\right) \ast G_{\nu,q}\\
\\
+\sum_{j=1}^D \left( \left( v^{f,-}_j\ast_{sp} G_{\nu}\right) \left(  v^{f,-}_i\ast_{sp}G_{\nu,j,p}\right) 
\right) \ast G_{\nu,q}
\end{array}
\end{equation}
Since we assume $\nu >0$, we have to show that for all $1\leq i,j\leq D$ the functions
\begin{equation}
{\Big |}B^0_{i,p,q}(t,x){\big |}\leq cr^{2\beta_0-3},
\end{equation}
and we can argue similarly as in the case of the Leray projection term in order to conclude that for a $t>0$
\begin{equation}
B^*_i(t,.)\in H^2\cap C^2.
\end{equation}
Moreover,as $\nu >0$ the nonlinear terms $B^*_i(t,.)$ and $L^*_i(t,.)$ and their derivatives up to second order are Lipschitz continuous (note the difference to the initial data which are only locally Lipschiitz continuous). It follows that after the first iteration the nonlinear terms $B^1_i(t,.)$ and $L^1_i(t,.)$ and their derivatives up to second order are Lipschitz continuous.
For $k=2$ (\ref{Navlerayscheme}) becomes
\begin{equation}\label{Navlerayscheme2}
\begin{array}{ll}
 v^{\nu,-,2}_i=v^{f,-}_i(.)\ast_{sp}G_{\nu}
+\sum_{j=1}^D v^{\nu,-,1}_j\frac{\partial v^{\nu,-,1}_i}{\partial x_j}\ast G_{\nu}\\
\\-\sum_{j,m=1}^D\int_{{\mathbb R}^D}\left( \frac{\partial}{\partial x_i}K_D(.-y)\right) \sum_{j,m=1}^D\left( \frac{\partial v^{\nu,-,1}_m}{\partial x_j}\frac{\partial v^{\nu,-,1}_j}{\partial x_m}\right) (t,y)dy\ast G_{\nu}\\
\\
=:v^{f,-}_i(.)\ast_{sp}G_{\nu}
+B_1\ast G_{\nu}-L_1\ast G_{\nu},
\end{array}
\end{equation}
where $B_1$ and $L_1$ denote abbreviations of the next order of approximation of the Burgers term and the Leray projection term. As we have $\left( v^{f,-}_i(.)\ast_{sp}G_{\nu}\right)(t,.)\in H^3 \cap C^3$ for $t>0$. and $B_1$ and $L_1$ are $C^2\cap H^2$, and  Lipschitz continuity of $B_1$ and $L_1$ and their derivatives up to second order we conclude that
\begin{equation}
\forall t>0~~ \left( B_1\ast G_{\nu}\right)(t,.)\in H^3 \cap C^3,~ \left( L_1\ast G_{\nu}\right)(t,.)\in H^3 \cap C^3.
\end{equation}

ad iii) For $k\geq 1$ choose a number $m$ is such that for $t\in [0,T]$
\begin{equation}
\forall~0\leq |\gamma|\leq m~D^{\gamma}_x\delta v^{init,\nu,-,k}_i(t,.) \mbox{is continuous and bounded.}
\end{equation}
We have observed that for $k=2$ we can choose $m=3$ such that the upper bound is independent of $\nu >0$.
For $0\leq |\gamma|\leq m$ and for $|\beta|+1=|\gamma|,~\beta_j+1=\gamma_j$ (if $|\gamma|\geq 1$) we consider a representation of  $D^{\gamma}_xv^{\nu,-,k}_i,~1\leq i\leq D,~k\geq 1$ of the form
\begin{equation}\label{Navlerayschemeiii}
\begin{array}{ll}
 D^{\gamma}_xv^{\nu,-,k}_i=v^{f,-}_i\ast_{sp}D^{\gamma}_xG_{\nu}
-D^{\beta}_x\left( \sum_{j=1}^D v^{\nu,-,k-1}_j\frac{\partial v^{\nu,-,k-1}_i}{\partial x_j}\right) \ast G_{\nu,j}\\
\\-D^{\beta}_x\left( \sum_{j,m=1}^D\int_{{\mathbb R}^D}\left( \frac{\partial}{\partial x_i}K_D(.-y)\right) \sum_{l,m=1}^D\left( \frac{\partial v^{\nu,-,k-1}_m}{\partial x_l}\frac{\partial v^{\nu,-,k-1}_l}{\partial x_m}\right) (t,y)dy\right) \ast G_{\nu,j}.
\end{array}
\end{equation}
Here, recall $G_{\nu}$ is the fundamental solution of the heat equation $p_{,t}-\nu\Delta p=0$, $\ast$ denotes the convolution, $\ast_{sp}$ denotes the spatial convolution, and $K_D$ denotes the fundamental solution of the 
Laplacian equation for dimension $D\geq 3$. In the following the constant $c>0$ is generic.
Note that for $1\leq i\leq D$ the initial data $v^{f,-}_i$have polynomial decay of any order at spatial infinity, i.e. we have for $|x|\geq 1$
\begin{equation}
{\big |}v^{f,-}_i(x){\big |}\leq \frac{c}{1+|x|^{2(D+1)+D+m}}
\end{equation}
 for some finite constant $c>0$ and $t\geq 0$. Hence, for multiindices $0\leq |\gamma|\leq m$ and $t>0$ we have for some finite constant $c>0$  and for for $|x|\geq 1$
\begin{equation}\label{I}
{\big |}v^{f,-}_i\ast_{sp}D^{\gamma}_xG_{\nu}(t,x){\big |}\leq \frac{c}{1+|x|^{2(D+1)}}.
\end{equation}
Assuming inductively 
\begin{equation}\label{indhyp}
\forall~l\leq k-1~\forall 0\leq|\gamma|\leq m~{\Big |}D^{\gamma}_xv^{\nu,-,l}_i(.){\Big |}\leq  \frac{c}{1+|x|^{2(D+1)}}
\end{equation}
we have or some finite constant $c>0$  and for for $|x|\geq 1$ 
\begin{equation}
{\big |}D^{\beta}_xB^{k-1}{\big |}:={\Big |}\sum_{j=1}^D D^{\beta}_x\left( v^{\nu,-,k-1}_j\frac{\partial v^{\nu,-,k-1}_i}{\partial x_j}\right) (t,.){\Big |}\leq \frac{c}{1+|x|^{4(D+1)}},
\end{equation}
where $|\beta|+1+|\gamma|$ is defined as above, and 
\begin{equation}
{\big |}D^{\beta}_xL^{k-1}{\big |}\leq \frac{c}{1+|x|^{4D+3}},
\end{equation}
where 
\begin{equation}
D^{\beta}_xL^{k-1}\equiv\sum_{j,m=1}^D\int_{{\mathbb R}^D}\left( \frac{\partial}{\partial x_i}K_D(.-y)\right) \sum_{j,m=1}^D\left( D^{\beta}_x\left( \frac{\partial v^{\nu,-,k-1}_m}{\partial x_j}\frac{\partial v^{\nu,-,k-1}_j}{\partial x_m}\right)\right)  (t,y)dy.
\end{equation}
Convolutions with  $G_{\nu}$ or $G_{\nu,i}$ weaken this polynomial decay by order $D$ at most such that
we (generously) get 
\begin{equation}\label{B}
{\big |}D^{\beta}_xB^{k-1}\ast G_{\nu,j}{\big |}\leq \frac{c}{1+|x|^{3(D+1)}}
\end{equation}
and
\begin{equation}\label{L}
{\Big |}D^{\beta}_xL^{k-1}\ast G_{\nu,j}{\Big |}\leq \frac{c}{1+|x|^{2(D+1)}}.
\end{equation}
Hence using the representation (\ref{Navlerayscheme}) and (\ref{I}),(\ref{B}),(\ref{L}) we get
\begin{equation}\label{ind}
\forall~l\leq k~\forall 0\leq|\gamma|\leq m~{\big |}D^{\gamma}_xv^{\nu,-,l}_i(.){\big |}\leq  \frac{c}{1+|x|^{2(D+1)}}
\end{equation}
and by (\ref{indhyp}) the same holds for the increments $D^{\gamma}_x\delta v^{\nu,-,k}_i(.)$ a fortiori.

Note that similar conclusions can be made using a vorticity iteration scheme $\left( \omega^{\nu,-,k}_i\right)_{k\geq 0,~1\leq i\leq D}$, where for $k=0$
\begin{equation}
 \omega^{\nu,-,0}_i=\omega^{f,-}_i,~\mbox{for}~1\leq i\leq D,
\end{equation}
and where for $k\geq 1$ the function $\omega^{\nu,-,k}_i,~1\leq i\leq D$ is determined as the time-local solution of the Cauchy problem 
\begin{equation}\label{vort2a}
\left\lbrace \begin{array}{ll}
 \frac{\partial \omega^{\nu,-,k}_i}{\partial s}-\nu\Delta \omega^{\nu,-,k}_i-\sum_{j=1}^3v^{\nu,-,k-1}_j\frac{\partial \omega^{\nu,-,k-1}_i}{\partial x_j}\\
 \\
 =-\sum_{j=1}^3\frac{1}{2}\left(\frac{\partial v^{\nu,-,k-1}_i}{\partial x_j}+\frac{\partial v^{\nu,-,k-1}_j}{\partial x_i} \right)\omega^{\nu,-,k-1}_j,\\
 \\
 \omega^{\nu,-,k}_i(0,.)=\omega^{f,-}_i,
 \end{array}\right.
\end{equation}
on a domain $[0,T]\times {\mathbb R}^D$ for some time horizon $T>0$.
Recall that in the limit $k\uparrow \infty$ (if existent) we have
\begin{equation}\label{veli0rk}
v^{\nu,-}(t,x)=\int_{{\mathbb R}^3}\left( K_3(y)\omega^{\nu,-}(t,x-y)\right) dy,
\end{equation}
by the Bio-Savart law .

ad iv) 
For $k\geq 3$ we have $\delta v^{\nu,-,k}_i(t,.)\in H^3\cap C^3$ corresponding to $\delta \omega^{\nu,-,k}_i(t,.)\in H^2\cap C^2$. For convenience of the reader we do some explicit very simple computations concerning the recursive relation of increments.
For the scheme $\omega^{\nu,-,k}_i(s,.),~k\geq 1$ with
\begin{equation}\label{vort2ahelp}
\begin{array}{ll}
   \omega^{\nu,-,k}_i(s,.)=\omega^{f,-}_i\ast_{sp}G_{\nu}+\left( \sum_{j=1}^3v^{\nu,-,k-1}_j\frac{\partial \omega^{\nu,-,k-1}_i}{\partial x_j}\right) \ast G_{\nu}\\
 \\
 -\left( \sum_{j=1}^3\frac{1}{2}\left(\frac{\partial v^{\nu,-,k-1}_i}{\partial x_j}+\frac{\partial v^{\nu,-,k-1}_j}{\partial x_i} \right)\omega^{\nu,-,k-1}_j\right)\ast G_{\nu},
 \end{array}
\end{equation}
we have
\begin{equation}\label{vort2ahelp3}
\begin{array}{ll}
   \omega^{\nu,-,k}_i(s,.)- \omega^{\nu,-,k-1}_i(s,.)=\left( \sum_{j=1}^3v^{\nu,-,k-1}_j\frac{\partial \omega^{\nu,-,k-1}_i}{\partial x_j}\right) \ast G_{\nu}\\
 \\
 -\left( \sum_{j=1}^3\frac{1}{2}\left(\frac{\partial v^{\nu,-,k-1}_i}{\partial x_j}+\frac{\partial v^{\nu,-,k-1}_j}{\partial x_i} \right)\omega^{\nu,-,k-1}_j\right)\ast G_{\nu}\\
 \\
-\left( \sum_{j=1}^3v^{\nu,-,k-2}_j\frac{\partial \omega^{\nu,-,k-1}_i}{\partial x_j}\right) \ast G_{\nu}\\
 \\
 +\left( \sum_{j=1}^3\frac{1}{2}\left(\frac{\partial v^{\nu,-,k-2}_i}{\partial x_j}+\frac{\partial v^{\nu,-,k-2}_j}{\partial x_i} \right)\omega^{\nu,-,k-1}_j\right)\ast G_{\nu} \\
 \\
 +\left( \sum_{j=1}^3v^{\nu,-,k-2}_j\frac{\partial \omega^{\nu,-,k-1}_i}{\partial x_j}\right) \ast G_{\nu}\\
 \\
 -\left( \sum_{j=1}^3\frac{1}{2}\left(\frac{\partial v^{\nu,-,k-2}_i}{\partial x_j}+\frac{\partial v^{\nu,-,k-2}_j}{\partial x_i} \right)\omega^{\nu,-,k-1}_j\right)\ast G_{\nu}
 \\
\\
-\left( \sum_{j=1}^3v^{\nu,-,k-2}_j\frac{\partial \omega^{\nu,-,k-2}_i}{\partial x_j}\right) \ast G_{\nu}\\
 \\
 +\left( \sum_{j=1}^3\frac{1}{2}\left(\frac{\partial v^{\nu,-,k-2}_i}{\partial x_j}+\frac{\partial v^{\nu,-,k-2}_j}{\partial x_i} \right)\omega^{\nu,-,k-2}_j\right)\ast G_{\nu}\\
\end{array}
\end{equation}
$$
\begin{array}{ll}
 =\left( \sum_{j=1}^3\left( v^{\nu,-,k-1}_j- v^{\nu,-,k-2}_j\right) \frac{\partial \omega^{\nu,-,k-1}_i}{\partial x_j}\right) \ast G_{\nu}\\
 \\
  +\sum_{j=1}^3v^{\nu,-,k-2}_j
  \left( \omega^{\nu,-,k-1}_j-\omega^{\nu,-,k-2}_j\right)\ast G_{\nu}\\
 \\
 -\left( \sum_{j=1}^3\left( \frac{1}{2}\left(\frac{\partial \delta v^{\nu,-,k-1}_i}{\partial x_j}+\frac{\partial \delta v^{\nu,-,k-1}_j}{\partial x_i} \right)\right) \omega^{\nu,-,k-1}_j\right)\ast G_{\nu}\\
 \\
 -\left( \sum_{j=1}^3\frac{1}{2}\left(\frac{\partial v^{\nu,-,k-2}_i}{\partial x_j}+\frac{\partial v^{\nu,-,k-2}_j}{\partial x_i} \right)\delta \omega^{\nu,-,k-1}_j\right)\ast G_{\nu}.
 \end{array}
$$
Hence the increments 
\begin{equation}
\delta \omega^{\nu,-,k}_i:=\omega^{\nu,-,k}_i-\omega^{\nu,-,k-1}_i
\end{equation}
satisfy the recursion
\begin{equation}\label{omegain}
\begin{array}{ll}
\delta \omega^{\nu,-,k}_i=&\left( \sum_{j=1}^3\delta v^{\nu,-,k-1}_j\frac{\partial \omega^{\nu,-,k-1}_i}{\partial x_j}\right)\ast G_{\nu}\\
\\
&+\left( \sum_{j=1}^3v^{\nu,-,k-1}_j\frac{\partial \delta \omega^{\nu,-,k-1}_i}{\partial x_j}\right)\ast G_{\nu}\\
\\
&-\left(  \sum_{j=1}^3\frac{1}{2}\left(\frac{\partial \delta v^{\nu,-,k-1}_i}{\partial x_j}+\frac{\partial \delta v^{\nu,-,k-1}_j}{\partial x_i} \right)\omega^{\nu,-,k-1}_j\right)\ast G_{\nu}\\
\\
&-\left(  \sum_{j=1}^3\frac{1}{2}\left(\frac{\partial  v^{\nu,-,k-2}_i}{\partial x_j}+\frac{\partial  v^{\nu,-,k-2}_j}{\partial x_i} \right)\delta \omega^{\nu,-,k-1}_j\right)\ast G_{\nu},
\end{array}
\end{equation}
where
\begin{equation}
\delta v^{\nu,-,k-1}_i:=v^{\nu,-,k-1}_i-v^{\nu,-,k-2}_i.
\end{equation}
Elementary arguments (we considered similar arguments elsewhere) lead to the first part of the next lemma. The second part ($\nu$-independent estimates) is proved in the next item.  
\begin{lem}
Given $\nu >0$ there exists $T>0$ and a constant $c\in (0,1)$ (dependent on $\nu$) such that for all $k\geq 3$ we have
\begin{equation}\label{vortcontr1}
\max_{1\leq i\leq D}\sup_{s\in [0,T]}{\big |}\delta \omega^{\nu,-,k+1}_i(s,.){\big |}_{H^2\cap C^2}\leq c\max_{1\leq i\leq D}\sup_{s\in [0,T]}{\big |}\delta \omega^{\nu,-,k}_i(s,.){\big |}_{H^2\cap C^2},
\end{equation}
and for all $k\geq 3$
\begin{equation}\label{velcontr1}
\max_{1\leq i\leq D}\sup_{s\in [0,T]}{\big |}\delta v^{\nu,-,k+1}_i(s,.){\big |}_{H^3\cap C^3}\leq c\max_{1\leq i\leq D}\sup_{s\in [0,T]}{\big |}\delta v^{\nu,-,k}_i(s,.){\big |}_{H^3\cap C^3}.
\end{equation}
As $\beta_0<2+\alpha_0$ is close to $2+\alpha_0$ these estimates have the $\nu$-independent counterparts with a loss of regularity.
There exists $T>0$ and a constant $c\in (0,1)$ independent on $\nu$ such that for all $k\geq 3$ we have
\begin{equation}\label{vortcontr2}
\max_{1\leq i\leq D}\sup_{s\in [0,T]}{\big |}\delta \omega^{\nu,-,k+1}_i(s,.){\big |}_{H^1\cap C^1}\leq c\max_{1\leq i\leq D}\sup_{s\in [0,T]}{\big |}\delta \omega^{\nu,-,k}_i(s,.){\big |}_{H^1\cap C^1},
\end{equation}
and for all $k\geq 3$
\begin{equation}\label{velcontr2}
\max_{1\leq i\leq D}\sup_{s\in [0,T]}{\big |}\delta v^{\nu,-,k+1}_i(s,.){\big |}_{H^2\cap C^2}\leq c\max_{1\leq i\leq D}\sup_{s\in [0,T]}{\big |}\delta v^{\nu,-,k}_i(s,.){\big |}_{H^2\cap C^2}.
\end{equation}
Moreover, the statement of (\ref{vortcontr1}) and (\ref{velcontr1}) hold for each $t_0>0$ in a time interval $[t_0,T]$ for some $T>0$.
\end{lem}

ad v) We have to supplement the previous argumentsin order to show  independence of $\nu>0$, and we have to show how a viscosity limit can be obtained which preserves essential properties such as contraction with lower regularity or contraction with higher regularity after finite time (as indicated in the last lemma). We have to show va) the estimate above are essentially independent of the viscosity $\nu >0$, and vb) that the viscosity limit preserves incompressibility, or that incompressibility is ontained in the viscosity limit. Note that some of the following estimates can be simplified in case of the data choice in item ib) where we have chosen the data $v^{f,-}_j,~j\neq i_0$ that are regular.

We note that in the case of the data choice in subitem ib) in item i) above we have regularity of the data $v^{f,-}_j,~j\neq i_0$, where $v^{f,-}_{i_0}$ is a H\"older continuous data function with singular 
vorticity. Especially we have $v^{f,-}_j\in H^1\cap C^1$for $j\neq i_0$. In this case we observe from  (\ref{Navlerayschemeiii}) thta the iteration scheme for the $D^{\gamma}_x\delta v^{\nu,-,k}_i$ increment 
bcomes for $k\geq 1$ 
\begin{equation}\label{Navlerayschemeiii}
\begin{array}{ll}
D^{\gamma}_x\delta v^{\nu,-,k}_i= D^{\gamma}_xv^{\nu,-,k}_i-v^{f,-}_i\ast_{sp}D^{\gamma}_xG_{\nu}\\
\\
=-D^{\beta}_x\left( \sum_{j=1}^D v^{\nu,-,k-1}_j\frac{\partial v^{\nu,-,k-1}_i}{\partial x_j}\right) \ast G_{\nu,j}-\\
\\D^{\beta}_x{\Big(} \sum_{j,m=1}^D\int_{{\mathbb R}^D}\left( \frac{\partial}{\partial x_i}K_D(.-y)\right)\sum_{l,m=1}^D\left( \frac{\partial v^{\nu,-,k-1}_m}{\partial x_l}\frac{\partial v^{\nu,-,k-1}_l}{\partial x_m}\right) (t,y)dy{\Big )} \ast G_{\nu,j}.
\end{array}
\end{equation}
For $|\gamma|\leq 1$ and $\beta=0$ this iteration scheme can be applied directly, where at the first step we apply convoluted data functions $v^{f,-}_{j}\ast_{sp} G_{\nu},~ 1\leq j\leq 3$ and  $v^{f,-}_{j,k}\ast_{sp} G_{\nu},~1\leq j,k\leq 3,~j\neq i_0$ with the substitution $v^{f,-}_{i_0,i_0}\ast G_{\nu}=-\sum_{j\neq i_0} v^{f,-}_{j,j}\ast G_{\nu}$ for $j\neq i_0$ and $t>0$. We then get a simple viscosity limit for the local solution (cf. also below).

In the general case, i.e., in any case of data choice considered in item i) we can use Lipschitz continuous upper bounds of the Leray  projection terms and the Burgers term. Note again that spatial convolutions with first order spatial derivatives of the Gaussian degenerate, but that the time- and spatial Gauusian first order derivatives convolutions of Lipschitz continuous Burgers- and Leray projection terms contribute to the solution (this is easily demonstrated by considering simple examples).

\begin{equation}\label{Navlerayscheme2v}
\begin{array}{ll}
 v^{\nu,-,2}_i=v^{f,-}_i(.)\ast_{sp}G_{\nu}
+\sum_{j=1}^D v^{\nu,-,1}_j\frac{\partial v^{\nu,-,1}_i}{\partial x_j}\ast G_{\nu}\\
\\-\sum_{j,m=1}^D\int_{{\mathbb R}^D}\left( \frac{\partial}{\partial x_i}K_D(.-y)\right) \sum_{j,m=1}^D\left( \frac{\partial v^{\nu,-,1}_m}{\partial x_j}\frac{\partial v^{\nu,-,1}_j}{\partial x_m}\right) (t,y)dy\ast G_{\nu}\\
\\
=:v^{f,-}_i(.)\ast_{sp}G_{\nu}
+B_1\ast G_{\nu}-L_1\ast G_{\nu},
\end{array}
\end{equation}
Ad va)  The main observation here is the estimate of 
\begin{equation}\label{obs}
v^f_{i_0}\ast_{sp}G_{\nu,j}
\end{equation}
for given $1\leq j \leq D$. First note that for $r=|x|$ we have an upper bound
\begin{equation}\label{dergnu}
{\big |}G_{\nu,j}(t,x){\big |}\leq \frac{c}{(\nu t)^{\delta}r^{D+1-2\delta}},
\end{equation}
where integrability is given for $\delta\in (0.5,1)$ and $c$ is independent of $\nu$. In (\ref{dergnu}) independence of $c$ from $\nu$ follows from
\begin{equation}
\begin{array}{ll}
{\big |}G_{\nu,k}(t,z){\big |}={\Big |}\left( \frac{-2z_k}{4\rho \nu t}\right) \frac{1}{\sqrt{4\pi \nu t}^D}\exp\left(-\frac{|z|^2}{4\nu t} \right){\Big |}\\
\\
\leq {\Big |}\left( \frac{-2z_k}{4\rho \nu t}\right)\exp\left(-\frac{|z|^2}{8\nu t} \right) \frac{\sqrt{2}}{\sqrt{8\pi \nu t}^D}\exp\left(-\frac{|z|^2}{8\nu t} \right){\Big |}\\
\\
\leq {\Big |}\left( \frac{2|z|^2}{4\rho \nu t}\right)\exp\left(-\frac{|z|^2}{8\nu t} \right) \frac{1}{|z|}\frac{\sqrt{2}}{\sqrt{8\pi \nu t}^D}\exp\left(-\frac{|z|^2}{8\nu t} \right){\Big |}\\
\\
\leq {\Big |}c_0 \frac{1}{|z|}\frac{\sqrt{2}}{\sqrt{8\pi \nu t}^D}\exp\left(-\frac{|z|^2}{8\nu t} \right){\Big |}
\end{array}
\end{equation}
for some finite constant $c_0$ which is independent of $\nu$. Then the standard pointwise upper bound estimate for $G_{\nu}$ can be employed and multiplied with $\frac{1}{|z|}$ in order to obtain that $c$ in (\ref{dergnu}) is independent of $\nu$.
Next
\begin{equation}\label{aber}
\begin{array}{ll}
{\big |}v^{f,-}_{i_0}\ast_{sp}G_{\nu,j}(t,x){\big |}\leq \int(r_0)^{\beta_0}\frac{\tilde{c}}{(\nu t)^{\delta}||r-r_0´|^{D+1-2\delta}}r_0^{D-1}dr_0\\
\\
\leq c^*+\int|r_0|^{\beta_0-1}\frac{\tilde{c}^*}{(\nu t)^{\delta}|r-r_0|^{D-2\delta}}r_0^{D-1}dr_0\\
\\
=c^*+\int|r_0-r|^{\beta_0-1}\frac{\tilde{c}^*}{(\nu t)^{\delta}|r_0|^{D-2\delta}}r_0^{D-1}dr_0
\end{array}
\end{equation}
upon partial integration in a sufficiently large ball and with some finite constants $c^*$ and $\tilde{c}^*$ which do not depend on $\nu$. For each $\nu >0$ we may consider the right side of the integral in a ball of radius $\sqrt{\nu}$. Outside such a ball we have strong decay due to the Gaussian. Hence we get
for given $1\leq j \leq D$ and  $r=|x|$ 
\begin{equation}\label{dergnu2}
{\big |}v^{f,-}_{i_0}\ast_{sp}G_{\nu,j}(t,x){\big |}\leq C^*r^{\beta_0-1},
\end{equation}
for a $C^*$ independent of $\nu$. Similarly we find some finite constant $C_+$ which is independent of $\nu >0$ such that 
\begin{equation}\label{dergnu3}
\begin{array}{ll}
{\big |}v^{f,-}_{i_0}\ast_{sp}G_{\nu}(t,x){\big |}\leq C_+r^{\beta_0},~ {\big |}\omega^{f,-}_{i_0}\ast_{sp}G_{\nu}(t,x){\big |}\leq C_+r^{\beta_0-(1+\alpha_0)}\\
\\
{\big |}\omega^{f,-}_{i_0}\ast_{sp}G_{\nu,j}(t,x){\big |}\leq C_+r^{\beta_0-1-(1+\alpha_0)}
\end{array}
\end{equation}Using these estimates we can use local time solution representations for the velocity and prove $v_i(t,.)\in H^2\cap C^2$ for $1\leq i\leq D$, or local time solution representations of the vorticity in order to prove $\omega_i(t,.)\in ,~ 1\leq i\leq D$ on some time interval $[t_0,t_0+\Delta]$. Note thatfor all local solution on a time interval $[t_0,t_0+\Delta]$ we have
\begin{equation}
v^{-}_i(t,.)\in H^2\cap C^2\Rightarrow \omega^{-}_i(t,.)\in H^1\cap C^1,
\end{equation}
and, 
\begin{equation}
\omega^{-}_i(t,.)\in H^1\cap C^1\Rightarrow v^{-}_i(t,.)\in H^2\cap C^2~\mbox{using Biot Savart.}
\end{equation}
We prefer to argue with the vorticity, but also mention briefly an alternative for the velocity. 
A local solution of the equation
\begin{equation}\label{vort2}
\begin{array}{ll}
 \frac{\partial \omega^{\nu,-}_i}{\partial \tau}
 -\nu \Delta \omega^{\nu,-}_i-\sum_{j=1}^3v^{\nu,-}_j\frac{\partial \omega^{\nu,-}_i}{\partial y_j}
=-\sum_{j=1}^3\frac{1}{2}\left(\frac{\partial v^{\nu,-}_i}{\partial x_j}+\frac{\partial v^{\nu,-}_j}{\partial x_i} \right)\omega^{\nu,-}_j
\end{array}
\end{equation}
may be represented in the form
\begin{equation}\label{vort2}
\begin{array}{ll}
 \omega^{\nu,-}_i
 = \omega^{f,-}_i\ast_{sp}G_{\nu}+\left( \sum_{j=1}^3v^{\nu,-}_j\frac{\partial \omega^{\nu,-}_i}{\partial y_j}\right) \ast G_{\nu}\\
 \\
-\left( \sum_{j=1}^3\frac{1}{2}\left(\frac{\partial v^{\nu,-}_i}{\partial x_j}+\frac{\partial v^{\nu,-}_j}{\partial x_i} \right)\omega^{\nu,-}_j\right) \ast G_{\nu}\\
\\
=\omega^{f,-}_i\ast_{sp}G_{\nu}-\left( \sum_{j=1}^3v^{\nu,-,}_{j,j}\omega^{\nu}_i\right) \ast G_{\nu}+\sum_{j=1}^3\left( v^{\nu,-}_{j}\omega^{\nu,-}_i\right) \ast G_{\nu,j}\\
\\
-\left( \sum_{j=1}^3\frac{1}{2}\left(\frac{\partial v^{\nu,-}_i}{\partial x_j}+\frac{\partial v^{\nu,-}_j}{\partial x_i} \right)\omega^{\nu,-}_j\right) \ast G_{\nu}\\
\\
=\sum_{j=1}^3\left( v^{\nu,-}_{j}\omega^{\nu,-}_i\right) \ast G_{\nu,j}-\left( \sum_{j=1}^3\frac{1}{2}\left(\frac{\partial v^{\nu,-}_i}{\partial x_j}+\frac{\partial v^{\nu,-}_j}{\partial x_i} \right)\omega^{\nu,-}_j\right) \ast G_{\nu},
\end{array}
\end{equation}
where we use incompressibility in the latter step. Next in our iteration scheme the first order approximation of the short time vorticity increment value is
\begin{equation}\label{vort2}
\begin{array}{ll}
 \delta \omega^{\nu,-,1}_i=\sum_{j=1}^3\left(  v^{f,-}_{j}\ast_{sp} G_{\nu}\right)\left(  \omega^{f,-}_i \ast_{sp} G_{\nu}\right) \ast G_{\nu,j}\\
 \\
 -\left(  \sum_{j=1}^3\frac{1}{2}\left(  v^{f,-}_i\ast_{sp}G_{\nu,j}+v^{f}_j\ast_{sp}G_{\nu,i}\right) \omega^{f,-}_j\ast_{sp}G_{\nu}\right)  \ast G_{\nu} ,
\end{array}
\end{equation}
For the first derivative of the increment we have
\begin{equation}\label{vort2}
\begin{array}{ll}
 \delta\omega^{\nu,-,1}_{i,k}= \sum_{j=1}^3  \left( \left(  v^{f,-}_{j}\ast_{sp} G_{\nu}\right)\left(  \omega^{f,-}_i \ast_{sp} G_{\nu}\right)\right)_{,k} \ast G_{\nu,j}\\
 \\
 -\left(  \sum_{j=1}^3\frac{1}{2}\left(  v^{f,-}_i\ast_{sp}G_{\nu,j}+v^{f,-}_j\ast_{sp}G_{\nu,i}\right) \omega^{f,-}_j\ast_{sp}G_{\nu}\right)  \ast G_{\nu,k} ,
\end{array}
\end{equation}
From the estimates in (\ref{dergnu2}) and (\ref{dergnu3}) it follows that for $\beta_0\in (2,2+\alpha_0)$ close to $2+\alpha_0$ the functions 
\begin{equation}
\begin{array}{ll}
\sum_{j=1}^3  \left( \left(  v^{f,-}_{j}\ast_{sp} G_{\nu}\right)\left(  \omega^{f,-}_i \ast_{sp} G_{\nu}\right)\right)_{,k}\\
\\
=\sum_{j=1}^3   \left(  v^{f,-}_{j}\ast_{sp} G_{\nu,k}\right)\left(  \omega^{f,-}_i \ast_{sp} G_{\nu}\right)\\
\\
+\sum_{j=1}^3   \left(  v^{f,-}_{j}\ast_{sp} G_{\nu}\right)\left(  \omega^{f,-}_i \ast_{sp} G_{\nu,k}\right)
\end{array}
\end{equation}
 are spatially Lipschitz such that
\begin{equation} 
\delta\omega^{\nu,-,1}_{i}(t,.)\in H^1\cap C^1.
\end{equation}
The contraction results in (\ref{vortcontr2}) and (\ref{velcontr2}) follow straightforwardly. For $t_0>0$ we can use this regularity at $t_0$ and may use the same iteration scheme starting in the interval $[t_0,t_0+\Delta_0]$ in order to obtain the contraction results in (\ref{vortcontr2}) and (\ref{velcontr2}) for some $T$ and all times $t\in [t_0,t_0+\Delta]$ for given $t_0 >0$and some $\Delta >0$.
Concerning the alternative argument for a local velocity solution, we first mention that the latter function satisfies
\begin{equation}\label{viii}
v_i^{\nu,-}(t,.)=v^{f,-}_i\ast G_{\nu}+\sum_{j=1}^Dv^{\nu,-}_jv^{\nu,-}_{i,j}\ast G_{\nu}+p^{\nu,-}_{,i}\ast G_{\nu}
\end{equation}
where $p^{\nu,-}_{,i}$ denotes the derivative of the pressure with respect to the variable $x_i$.
Now in a first approximation of our iteration scheme we have
\begin{equation}
\Delta p^{\nu,-,1}=\sum_{j,k=1}^D\left( v^{f,-}_{j}\ast_{sp}G_{\nu,k}\right) \left( v^{f,-}_{k}\ast_{sp} G_{\nu,j}\right) ,
\end{equation}
where the right side is Lipschitz as a production of two Lipschitz functions according to the estimates above. Using symmetric data we can conclude that for all $1\leq i\leq D$ the second derivatives  $p^1_{,i,i}$ of the first approximation of the pressure are Lipschitz. This leads to the conclusion that
\begin{equation}
p^1_{,i,i}\ast G_{\nu ,k}(t,.)\in L^2\cap C
\end{equation}
for some second order derivatives. This information can be used together with the  estimates in (\ref{dergnu2}) and (\ref{dergnu3}) to conclude that $v_i^{\nu_k,-}(t,.)\in C^2\cap H^2$ independently of $\nu_k$ for a sequence $(\nu_k)_k$ with $\nu_{k}\downarrow 0$.
Ad vb) for $1\leq i\leq D$ consider the local solution of the Navier Stokes type equation (related to the time reversed Euler equation) 
\begin{equation}\label{Navlerayschemeloc}
\begin{array}{ll}
 v^{\nu,-}_i=v^{\nu,-}_i(.)\ast_{sp}G_{\nu}
-\sum_{j=1}^D v^{\nu,-}_j\frac{\partial v^{\nu,-}_i}{\partial x_j}\ast G_{\nu}\\
\\-\sum_{j,m=1}^D\int_{{\mathbb R}^D}\left( \frac{\partial}{\partial x_i}K_D(.-y)\right) \sum_{j,m=1}^D\left( \frac{\partial v^{\nu,-}_m}{\partial x_j}\frac{\partial v^{\nu,-}_j}{\partial x_m}\right) (t,y)dy\ast G_{\nu}
\end{array}
\end{equation}
on the time integral $[0,T]$. We have to show that the incompressibility condition holds on the restricted time interval $(0,T]$, i.e.,
\begin{equation}
\mbox{for all $t\in (0,T]$}~\sum_{i=1}^Dv^{\nu,-}_{i,i}(t,.)=0
\end{equation}
The local solution in (\ref{Navlerayschemeloc}) satisfies the equation
 \begin{equation}\label{Navlerayeqincompressvb}
\begin{array}{ll}
 \frac{\partial v^{\nu,-}_i}{\partial t}-\nu\Delta v^{\nu,-}_i=
-\sum_{j=1}^D v^{\nu,-}_j\frac{\partial v^{\nu,-}_i}{\partial x_j}\\
\\+\sum_{j,m=1}^D\int_{{\mathbb R}^D}\left( \frac{\partial}{\partial x_i}K_D(.-y)\right) \sum_{j,m=1}^D\left( \frac{\partial v^{\nu,-}_m}{\partial x_j}\frac{\partial v^{\nu,-}_j}{\partial x_m}\right) (t,y)dy.
\end{array}
\end{equation} 
For $t\in (0,T]$ $v^{\nu,-}_i(t,.)\in H^3\cap C^3$ such that we can
apply the divergence operator in a classical sense to the equation in (\ref{Navlerayeqincompressvb}). We get
\begin{equation}\label{nsrev}
\begin{array}{ll}
 \frac{\partial \sum_{i=1}^Dv^{\nu,-}_{i,i}}{\partial t}-\nu \Delta \sum_{i=1}^D v^{\nu,-}_{i,i}
-\sum_{i=1}^D\sum_{j=1}^D v^{\nu,-}_{j,i}\frac{\partial \sum_{i=1}^Dv^{\nu,-}_{i,i}}{\partial x_j} \\
 \\
-\sum_{i=1}^D\sum_{j=1}^D v^{\nu,-}_{j,i}v^{\nu,-}_{i,j}\\
\\+\sum_{i=1}^D\sum_{j,m=1}^D\int_{{\mathbb R}^D}\left(K_{D,i,i}(.-y)\right) \sum_{j,m=1}^D\left( \frac{\partial v^{\nu,-}_m}{\partial x_j}\frac{\partial v^{\nu,-}_j}{\partial x_m}\right) (t,y)dy.
\end{array}
\end{equation}
The last two terms  in (\ref{nsrev}) encode the Poisson equation $\Delta p^-=\sum_{i=1}^D\sum_{j=1}^D v^{\nu,-}_{j,i}v^{\nu,-}_{i,j}$ (with $p^-$ the analogue of pressure), and we have $\sum_{i=1}^Dv^{f,-}_{i,i}=0$ in $H^1$ sense. Hence, $\sum_{i=1}^Dv^{\nu,-}_{i,i}(t,.)=0$ in classical sense is consistent with  (\ref{Navlerayeqincompressvb}).
Note that we have chosen $v^{f,-}_{j},~ 1\leq i\leq D,~ j\neq i_0$ at item i) above such that
\begin{equation}
\sum_{i=1}^Dv^{f,-}_i\ast_{sp} G_{\nu,i}(t,.)=0.
\end{equation}
We can use the local iteration scheme and check directly that the incompressibility condition holds, for the limit, i.e., $\sum_{i=1}^Dv^{\nu,-}_{i,i}(t,.)=0$.

ad vi) Choose a time horizon $T>0$ as in the previous step such that contraction holds for the higher order increments $\delta v^{\nu,-,k}_i$ with $k\geq 3$ as in (\ref{vel}). 
We choose the increment $\delta v^{\mbox{init}\nu,-,2}_i$ obtained after $k=2$ iterations of the scheme. We observed that there is an upper bound $C>0$ independent of $\nu$ such that
\begin{equation}
\begin{array}{ll}
{\big |}\delta v^{\mbox{init},\nu,-,2}_i+\sum_{l=3}^{\infty}\delta v^{\nu,-,l}_i{\big |}_{H^2\cap C^2}\\
\\
={\big |}v^{\nu,-}_i(t,.)-v^{f,-}_i\ast_{sp}G_{\nu}(t,.){\big |}_{H^2\cap C^2}\leq C
\end{array}
\end{equation}
where $C>0$ is independent of $\nu>0$. For all $\nu >0$ the function
\begin{equation}
v^{\nu,-}_i(t,.)=v^{f,-}_i\ast_{sp}G_{\nu}(t,.)+\delta v^{\mbox{init},\nu,-,2}_i+\sum_{l=3}^{\infty}\delta v^{\nu,-,l}_i
\end{equation}
satisfies the Navier stokes equation on the time interval $[0,T]$ such that
\begin{equation}\label{Navlerayeq2}
\begin{array}{ll}
\nu\Delta v^{\nu,-}_i= \frac{\partial v^{\nu,-}_i}{\partial t}
-\sum_{j=1}^D v^{\nu,-}_j\frac{\partial v^{\nu,-}_i}{\partial x_j}\\
\\+\sum_{j,m=1}^D\int_{{\mathbb R}^D}\left( \frac{\partial}{\partial x_i}K_D(.-y)\right) \sum_{j,m=1}^D\left( \frac{\partial v^{\nu,-}_m}{\partial x_j}\frac{\partial v^{\nu,-}_j}{\partial x_m}\right) (t,y)dy,
\end{array}
\end{equation}
where the right side of (\ref{Navlerayeq2}) is the time-reversed Euler equation operator applied to $v^{\nu,-}_i$. Recall that the time horizon $T>0$ is structurally independent of $\nu$ in the contraction result.
Next we may use the strong polynomial decay at spatial infinity in order to apply a strong compactness argument. Here we note that Rellich's theorem holds only with restrictions for unbounded domains. It is convenient to use the strong polynomial decay at spatial infinity in order to transform on a bounded domain. Moreover spatial transformations to bounded domains allow for convergence constructions in $C^m$-Banach spaces. Here recall that
\begin{prop}
For open and bounded $\Omega\subset {\mathbb R}^n$ and consider the function space
\begin{equation}
\begin{array}{ll}
C^m\left(\Omega\right):={\Big \{} f:\Omega \rightarrow {\mathbb R}|~\partial^{\alpha}f \mbox{ exists~for~}~|\alpha|\leq m\\
\\
\mbox{ and }\partial^{\alpha}f \mbox{ has an continuous extension to } \overline{\Omega}{\Big \}}
\end{array}
\end{equation}
where $\alpha=(\alpha_1,\cdots ,\alpha_n)$ denotes a multiindex and $\partial^{\alpha}$ denote partial derivatives with respect to this multiindex. Then the function space $C^m\left(\overline{\Omega}\right)$ with the norm
\begin{equation}
|f|_m:=|f|_{C^m\left(\overline{\Omega}\right) }:=\sum_{|\alpha|\leq m}{\big |}\partial^{\alpha}f{\big |}
\end{equation}
is a Banach space. Here,
\begin{equation}
{\big |}f{\big |}:=\sup_{x\in \Omega}|f(x)|.
\end{equation}

\end{prop}  
We choose a sequence $(\nu_p)_{p\geq 1}$ converging to zero and consider the spatial transformation
\begin{equation}
\delta v^{c,\mbox{init},-,\nu_p,2}_i(t,y)=\delta v^{\mbox{init},-,\nu_p,2}_i(t,x)
\end{equation}
for $y_j=\arctan(x_j),~1\leq j\leq D$ and for all $t\in [0,T]$. Note that for
\begin{equation}
{\big |}\delta v^{\mbox{init},-,\nu_p,2}_i(t,x){\big |}\leq \frac{c}{1+|x|^{2m}}
\end{equation}
with a finite $\nu$-independent constant $c$ and multiindices $\gamma$ with $0\leq |\gamma|\leq 3+\epsilon$ for some $\epsilon >0$ we have for all $t\in [0,T]$ and all $x\in {\mathbb R}^D$
\begin{equation}
{\big |}D^{\gamma}_y\delta v^{c,\mbox{init},-,\nu_p,2}_i(t,y){\big |}\leq c_0(1+|x|^{2m}){\big |}D^{\gamma}_x\delta v^{\mbox{init},-,\nu_p,2}_i(t,x){\big |}\leq C
\end{equation}
for some finite constants $c_0,C$. This implies 
\begin{equation}
\delta v^{c,\mbox{init},-,0,2}_i(t,.):=\lim_{\nu_p\downarrow 0}\delta v^{c,\mbox{init},-,\nu_p,2}_i(t,.)\in H^{2}\cap C^2~\mbox{ for all $0\leq t\leq T$},
\end{equation}
and as for some finite $C>0$ independent of $\nu_p$
\begin{equation}
\sup_{\nu_l>0}{\big |}(1+|x|^{2m}){\big |}D^{\gamma}_x\delta v^{\mbox{init},-,\nu_p,2}_i(t,x){\big |}\leq C
\end{equation}
we conclude that for all $t\in [0,T]$
\begin{equation}
\delta v^{c,\mbox{init},-,0,2}_i(t,.)\in H^2\cap C^2.
\end{equation}
The latter statement transfers to $\delta v^{\mbox{init},0,-,2}_i(t,.),~1\leq i\leq D$.
Similarly for the higher order increments $\delta v^{\nu,-,k}_i,~1\leq i\leq D$ for $k\geq 3$.
Hence the viscosity limit $v^{-}_i,~1\leq i\leq D$ satisfies 
\begin{equation}
v^{-}_i(t,.)-v^{f,-}_{i}\ast_{sp}G_{\nu}(t,.)=\delta v^{-}_i(t,.)\in H^2\cap C^2
\end{equation}
for all $t\in [0,T]$, where indeed
\begin{equation}
\delta v^{-}_i:=\lim_{\nu_p\downarrow 0 }\delta v^{\mbox{init},\nu_p,-,2}_i+\sum_{l=2}^{\infty}\delta v^{\nu_l,-,l}_i\in C^1\left(\left(0,T\right] , H^2\cap C^2\right). 
\end{equation} 
Next, we observe that the function $v^{-}_i,~1\leq i\leq D$ of the described regularity is a classical solution of the time-reversed Euler equation. Indeed, for some finite $C>0$
\begin{equation}
\lim_{\nu_p\downarrow 0}\nu_p{\big |}\Delta \delta v^{-}_i(t,.){\big |}_{C\cap L^2}\leq \lim_{\nu_p\downarrow 0}\nu_p\sup_{\nu_p}{\big |}\delta v^{\nu_p,-}_i(t,.){\big |}_{H^2\cap C^2}\leq \nu_pC\downarrow 0~\mbox{as $\nu_p\downarrow 0$}
\end{equation}
for the left side in (\ref{Navlerayeq2}), and hence the right side of (\ref{Navlerayeq2}) is also $0$.
Note that the set of continuous  functions on ${\mathbb R}^D$ which vanish at spatial infinite is closed. In our context of functions with very strong polynomial decay we may even transform to a bounded domain. Hence we may consider norms $|f|_{C^0}=\sup_x|f(x)|$ and similar norms $|f|_{C^m}$ for derivatives up to order $1\leq m\leq 3$ in our context.  
Furthermore, we have
\begin{equation}
\lim_{\nu_p\downarrow 0}\nu_p{\big |}\Delta (v^{f,-}_i\ast_{sp}G_{\nu})(t,.){\big |}_{C^0}
=\lim_{\nu_p\downarrow 0}\nu_p^{1-\delta}{\big |}\sum_{j} v^{f,-}_{i,j}\ast_{sp}\nu^{\delta}G_{\nu,j}(t,.){\big |}_{C^0}\downarrow 0 ~\mbox{as $\nu_p\downarrow 0$},
\end{equation}
where in the last step we may use the local upper bound
\begin{equation}
\begin{array}{ll}
|\int_{|x-y|\leq 1}(\sum_{j}v^{f,-}_{i,j})(t,y)\nu^{\delta}G_{\nu,j}(t,x;0,y)|\leq  \int c_0r^{\beta_0-(2+\alpha_0)}\frac{c_1}{t^{\delta}|x-y|^{4-2\delta}}dy\\
\\
\leq c_2r^{\beta_0+(2\delta-1) -(2+\alpha_0)}~ \mbox{for $0\leq r=\sqrt{y_1^2+y_2^2+y_3^2}$}
\end{array}
\end{equation}
in a ball $B_1(x)$ of radius $1$ around $x$, $\delta\in (0.5,1)$, and for some finite constants $c_0,c_1,c_2$ which are independent of $\nu$. For $\beta_0$ close to $2+\alpha_0$ this integral is bounded. Note that this is also achieved for weaker parameter conditions (for $\alpha_0,\beta_0$) if $\delta$ is chosen to be $1-\epsilon$ for small $\epsilon >0$. Note that the integral outside this ball clearly converges to zero as $\nu_p\downarrow 0$.  Here, note that we used the factor $\nu^{\delta}$  for the latter conclusion in order to ensure that a standard estimate of the Gaussian is independent of $\nu$.
Similar arguments show that the velocity viscosity limit $v^{-}$ has spatial regularity $H^3\cap C^3$ for positive time corresponding to a spatial vorticity regularity $H^2\cap C^2$ for positive time $t>0$.
We conclude that the original Euler equation develops in opposite time direction from data at time $t=0$ (corresponding to time $s=T$ of the reversed-time Euler equation) a weak singularity at time $T>0$ (corresponding to data at time $s=0$ the time-reversed Euler equation).


\begin{thebibliography}{19}
\baselineskip=12pt

%
%
 \bibitem{Ki}
 {\sc Kiselev, A.} {\em
 Blow-up for the 2D Euler equation on some bounded domains}, (to appear in Journal of Differential Equations).
\bibitem{LL}
{\sc Landau, L., Lifschitz, E.} {\em Lehrbuch der Theoretischen Physik VI, Hydrodynamik}, Akademie Verlag, Berlin. J., (1978).
\bibitem{MB}  
{\sc Majda, A., Bertozzi, L.}
{\em Vorticity and Incompressible Flow (Cambridge Texts in Applied Mathematics)}  Cambridge University Press , 2001.

%
%
 \end{thebibliography}
\end{document}